\documentclass[final,11pt]{elsarticle}
\graphicspath{{Fig/}}
\usepackage[utf8]{inputenc} 
\usepackage{graphicx} 
\usepackage{multirow}
\usepackage[inline]{changes}
\usepackage{booktabs} 
\usepackage{array} 
\usepackage{paralist} 
\usepackage{verbatim} 
\usepackage{tikz}
\usepackage{comment}
\usepackage{epstopdf}
\usepackage{epsfig}
\usepackage{diagbox}
\usepackage{mathtools}
\usepackage{algorithm}
\usepackage{algorithmicx}
\usepackage{algpseudocode}
\usetikzlibrary{arrows,shapes,chains}
\usepackage{amsmath,dsfont}
\usepackage{braket,amsfonts}
\usepackage{amsmath,dsfont}
\usepackage{harpoon}
\usepackage{amssymb}
\usepackage{mathrsfs}
\usepackage{graphics}
\usepackage{color}
\usepackage{geometry}
\usepackage{epsfig}
\usepackage{diagbox}
\usepackage{mathtools}
\usepackage{bbm}
\usepackage{xfrac}
\usepackage{fancyhdr}
\usepackage{graphicx}
\usepackage{wrapfig}
\usepackage{dsfont}
\usepackage{makecell}
\usepackage{tikz}
\usepackage{enumerate}
\usepackage{bm}
\usepackage{xcolor}
\usepackage{setspace}
\usepackage{array}

\def\a{\alpha}
\def\b{\beta}

\def\dd{\mbox{d}}

\def\x{{\bf x}}

\usepackage{amsthm}

\def\f{\frac}
\def\p{\partial}

\usepackage[caption=false]{subfig}
\usepackage{wrapfig}
\usepackage{dsfont}
\usepackage{enumerate}
\usepackage{xcolor}
\usepackage{setspace}
\usepackage{pgfplots}

\newtheorem{theorem}{Theorem}
\newtheorem{example}{\textup{\textbf{Example}}}

\usepackage{amssymb}

\journal{Applied Numerical Mathematics}

\begin{document}

\begin{frontmatter}



\title{Adaptive Hermite spectral methods in unbounded domains}


\author[UCLA]{Tom Chou}
\ead{tomchou@ucla.edu}

\author[PKU]{Sihong Shao}
\ead{sihong@math.pku.edu.cn}

\author[UCLA]{Mingtao Xia \corref{cor1}}
\ead{xiamingtao97@g.ucla.edu}

\cortext[cor1]{Corresponding author}

\address[UCLA]{Department of Mathematics, UCLA, Los Angeles, CA
  90095-1555, USA}
\address[PKU]{CAPT, LMAM and School of Mathematical
  Sciences, Peking University, Beijing 100871, CHINA}


\begin{abstract}
Recently, new adaptive techniques were developed that greatly improved
the efficiency of solving PDEs using spectral methods. These adaptive
spectral techniques are especially suited for accurately solving
problems in unbounded domains and require the monitoring and dynamic
adjustment of three key tunable parameters: the scaling factor, the
displacement of the basis functions, and the spectral expansion
order. There have been few analyses of numerical methods for unbounded
domain problems.  Specifically, there is no analysis of adaptive
spectral methods to provide insight into how to increase efficiency
and accuracy through dynamical adjustment of parameters. In this
paper, we perform the first numerical analysis of the adaptive
spectral method using generalized Hermite functions in both one- and
multi-dimensional problems. \added{Our analysis reveals why adaptive
spectral methods work well when a ``frequency indicator'' of the
numerical solution is controlled.  We then investigate how the
implementation of the adaptive spectral methods affects numerical
results, thereby providing guidelines for the proper tuning of
parameters.} Finally, we further improve performance by extending the
adaptive methods to allow bidirectional basis function translation,
and the prospect of carrying out similar numerical analysis to solving
PDEs arising from realistic difficult-to-solve unbounded models with
adaptive spectral methods is also briefly discussed.
\end{abstract}



\end{frontmatter}



\section{Introduction}
\label{sec:intro}

Unbounded domain problems require efficient numerical methods for
computation. For example, resolving the decay of the solution of
Schr\"odinger's equations at infinity requires efficient unbounded
domain algorithms \citep{li2018stability}. In population dynamics,
tracking cell volume blowup in structured population PDE models
demands high-accuracy numerical methods in unbounded domains
\citep{xia2021kinetic,xia2020pde}. Furthermore, in solid-state physics,
numerical methods for unbounded domains are required for studying
long-range particle interactions
\citep{hugli2012artificial,mengotti2011real}. Despite these numerous
applications, there has been little research on developing efficient
and accurate algorithms for solving models in unbounded domains.

Adaptive methods, such as re-defining grids for finite difference
methods \citep{ren2000iterative} and re-generating meshes for finite
element methods \citep{antonietti2019adaptive,
  babuska2012modeling,li2002adaptive,tang2003adaptive}, which are
applied to PDEs defined on finite domains, can dramatically improve
not only accuracy but computational efficiency.  Recently, novel
adaptive techniques for spectral methods have been developed and
incorporated into efficient algorithms for numerically solving PDEs in
unbounded domains that posed substantial numerical difficulties when
using previous numerical methods \citep{xia2020b,xia2020a}.
\added{The adaptive spectral methods consist of three separate but
  interdependent procedures: (i) a scaling technique that adjusts the
  shape of the basis functions to capture the varying decay rate of
  the function at infinity, (ii) a moving technique that adjusts the
  displacement of the basis function to better assign allocation
  points and capture intrinsic translation of the solution, and (iii)
  a $p$-adaptive technique that adjusts the expansion order of the
  numerical solution to deal with oscillations of the solution.}
These adaptive spectral techniques require tuning of three key
parameters: the scaling factor $\beta$, the displacement of the basis
function $x_0$, and the spectral expansion order $N$. For example, if
we use the generalized Hermite functions \citep{xiang2010generalized}
as basis functions on $\mathbb{R}$, the variables $\beta, x_0$, and
$N$ appear in a spectral expansion according to

\begin{equation}
U_{N, x_0}^{\beta}\coloneqq 
\sum_{i=0}^N u_{i, x_0}^{\beta}\hat{\mathcal{H}}_i^{\beta}(x-x_0) = 
\sum_{i=0}^N u_{i, x_0}^{\beta}\hat{\mathcal{H}}_i(\beta(x-x_0)),
\end{equation}
where $u_{i, x_0}^{\beta}$ is the coefficient of the $i^{\rm
  th}$-order generalized Hermite function $\hat{\mathcal{H}}_i^{\beta}$

\begin{equation}
  \hat{\mathcal{H}}_i^{\beta} \coloneqq \frac{1}{\sqrt{2^ii!}}
  H_i(\beta x)e^{-\tfrac{(\beta{x})^2}{2}},\,\,\,
H_i(x) = (-1)^i e^{x^2}\partial_x^i (e^{-x^2}).
\end{equation}  
For example, for PDEs involving a spatial variable $x\in\mathbb{R}$
and a temporal variable $t\in[0, T]$, we typically impose a spectral
expansion using generalized Hermite functions of $x$ and forward time
$t$ starting from an initial condition at $t=0$.

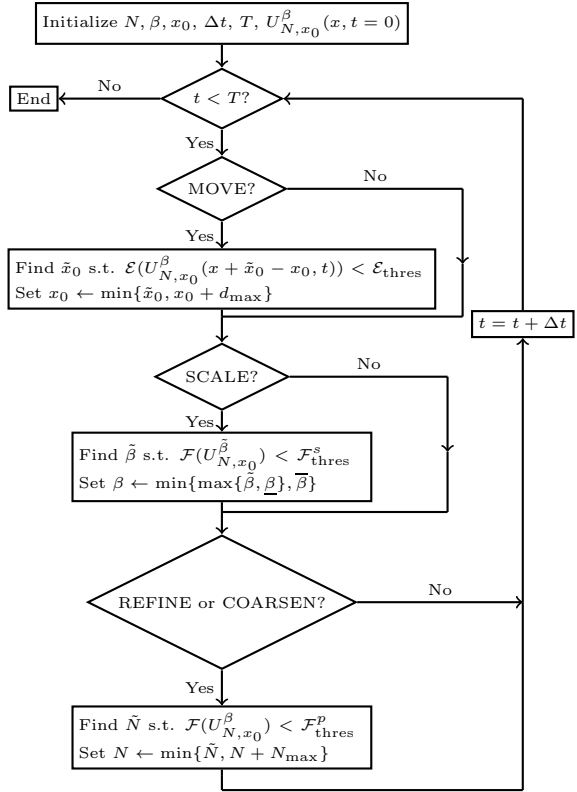
\begin{figure}  
\vspace{-0.15in}
\centering
		\tiny
		\tikzstyle{format}=[rectangle,draw,thick,fill=white]  
		\tikzstyle{test}=[diamond,aspect=2,draw,thick]  
		\tikzstyle{point}=[coordinate,on grid,thick]  
		\begin{tikzpicture}
		\node[format] (start){Initialize $N, \beta, x_0$, $\Delta t$, $T$, $U_{N, x_0}^{\beta}(x, t=0)$};
		\node[test,below of=start,node distance=10mm](time){$t<T$?};
		\node[format,left of=time,node distance=25mm](over){End};
		\node[test,below of=time,node distance=12mm] (PCFinit){{MOVE?}};
		\node[point, left of=PCFinit, node distance = 35mm](moveinter1){};
		\node[format,below of=PCFinit, text width  = 55mm, node distance=12mm](movee0){Find $\tilde{x}_0$ s.t. $\mathcal{E}(U_{N, x_0}^{\beta}(x+\tilde{x}_0-x_0, t))<\mathcal{E}_{\text{thres}}$ Set $x_0\gets \min\{\tilde{x}_0, x_0+d_{\max}\}$};
		\node[point, below of=movee0, node distance=5mm](pointinter1){};
		\node[point, right of=pointinter1, node distance=32mm](pointinter2){};
		\node[test,below of=pointinter1,node distance=8mm] (DS18init){SCALE?};
		\node[point, left of=DS18init, node distance=35mm](scaleinter){};
		\node[format,below of=DS18init, text width  = 38mm,node distance=12mm] (LCDinit){Find $\tilde{\beta}$ s.t. $\mathcal{F}(U_{N, x_0}^{\tilde{\beta}})<\mathcal{F}^s_{\text{thres}}$ Set $\beta\gets \min\{\max\{\tilde{\beta}, \underline{\beta}\}, \overline{\beta}\}$};
		\node[point, below of=LCDinit, node distance=6mm](pointscale){};
		\node[point, right of=pointscale, node distance = 30mm](pointscale2){};
		\node[test,below of=LCDinit,node distance=18mm](setkeycheck){REFINE or COARSEN?};
		\node[point,right of=PCFinit,node distance=32mm](movepoint){};
		\node[point,below of=movepoint,node distance=10mm](movepoint2){};
		\node[point,left of=movepoint2,node distance=18mm](movepoint3){};
		\node[point,left of=setkeycheck,node distance=15mm](point3){};
		\node[point,right of=point3,node distance=25mm](point4){};
		\node[point, right of=setkeycheck, node distance=40mm](pointconnect1){};
		\node[format, above of=pointconnect1, node distance=37mm](increasetime){$t = t + \Delta{t}$};
	     \node[point, above of=increasetime, node distance=30mm](pointconnect2){};
		\node[format,below of=setkeycheck, text width  = 38mm, node distance=18mm](processtime1){Find $\tilde{N}$ s.t. $\mathcal{F}(U_{N, x_0}^{\beta})<\mathcal{F}^p_{\text{thres}}$
 Set $N\gets \min\{\tilde{N}, N+N_{\max}\}$};
		\node[point, right of=DS18init, node distance=30mm](scalepoint1){};
		\node[point, below of=scalepoint1, node distance=10mm](scalepoint2){};
		\node[point, left of=scalepoint2, node distance=18mm](scalepoint3){};
		\node[point,below of=processtime1,node distance=5mm](display){};
		\node[point,below of=display,node distance=2mm](point1){};
		\node[point, right of=point1, node distance=40mm](refineover){};
		\draw[->, thick](start)--(time);
		\draw[-, thick](movee0)--(pointinter1);
		\draw[-, thick](LCDinit)--(pointscale);
		\draw[-, thick](pointinter1)--(pointinter2);
		\draw[->, thick](time)--node[left]{Yes}(PCFinit);
		\draw[->, thick](time)--node[above]{No}(over);
		\draw[->, thick](pointconnect1)--(increasetime);
		\draw[-, thick](increasetime)--(pointconnect2);
		\draw[-, thick](PCFinit)--node[above]{No}(movepoint);
		\draw[->, thick](movepoint)--(movepoint2);
		\draw[->, thick](DS18init)--node[left]{Yes}(LCDinit);
		\draw[->, thick](PCFinit)--node[left]{Yes}(movee0);
		\draw[->, thick](pointinter1)--(DS18init);
		\draw[-, thick](processtime1)--(display);
		\draw[-, thick](movepoint2)--(pointinter2);
		\draw[-, thick](DS18init)--node[above]{No}(scalepoint1);
		\draw[->, thick](scalepoint1)--(scalepoint2);
		\draw[-, thick](scalepoint2)--(pointscale2);
		\draw[->, thick](pointscale)->(setkeycheck);
		\draw[-, thick](pointscale2)--(pointscale);
		\draw[->, thick](setkeycheck)--node[above]{No}(pointconnect1);
		\draw[->, thick](pointconnect2)--(time);
		\draw[->, thick](setkeycheck) --node[left]{Yes}(processtime1);
		\draw[-, thick](display)--(point1);
		\draw[-, thick](point1)--(refineover);
		\draw[-, thick](refineover)--(pointconnect1);
		\end{tikzpicture}  
\caption{\footnotesize \added{Flow chart of an adaptive Hermite
    spectral method equipped with scaling, moving, and $p$-adaptive
    techniques. $x_0$ and $\tilde{x}_0$ are the displacements before
    and after the moving technique is used, and $\tilde{x}_0$ is
    chosen such that the exterior-error indicator
    $\mathcal{E}(U^{\beta}_{\tilde{x}_0, N})$ is below a moving
    threshold $\mathcal{E}_{\text{thres}}$. $\beta$ and
    $\tilde{\beta}$ are the scaling factors before and after scaling
    when the scaling technique is used, where the scaling factor
    $\tilde{\beta}$ is chosen such that
    $\mathcal{F}(U^{\tilde{\beta}}_{x_0, N})$ is below a scaling
    threshold $\mathcal{F}^s_{\text{thres}}$. $N$ and $\tilde{N}$ are
    the expansion orders before and after adjusting the expansion
    order when the $p$-adaptive technique is used. $\tilde{N}$ is
    chosen such that $\mathcal{F}(U^{\beta}_{x_0, \tilde{N}})$ is
    below a $p$-adaptivity threshold
    $\mathcal{F}^p_{\text{thres}}$. The three thresholds
    $\mathcal{E}_{\text{thres}}, \mathcal{F}^s_{\text{thres}}$, and
    $\mathcal{F}^p_{\text{thres}}$ are updated dynamically as time
    progresses. Details are described in \cite{xia2020b,xia2020a}. In
    addition, we impose constraints on the maximum allowable
    displacement $d_{\max}$ and expansion order increment $N_{\max}$
    within a single step, and the minimum and maximum scaling factors
    $\underline{\beta}$ and $\overline{\beta}$.}}
\label{algmovingscaling}
\vspace{-0.1in}
\end{figure}

\added{Adaptive spectral techniques are implemented as shown in
Fig.~\ref{algmovingscaling}. Specifically, the algorithm changes the
displacement of the basis function $x_0$ to control an
\textit{exterior-error indicator} $\mathcal{E}(U_{N, x_0}^{\beta})$
that reflects the ratio of the numerical solution's error outside a
given domain to the error in the whole domain. It also changes the
scaling factor $\beta$ as well as the spectral expansion order $N$ to
control a \textit{frequency indicator} $\mathcal{F}(U_{N,
  x_0}^{\beta})$ that measures the spread and oscillation of the
solution. The indicators are defined in \cite{xia2020b} as}

\begin{equation}
   \mathcal{E}(U_{N,x_0}^{\b}) = \f{\|\p_x U_{N,x_0}^{\b}
     \cdot\mathbb{I}_{(x_R, \infty)}\|}{\|\p_x U_{N,x_0}^{\b}\cdot
     \mathbb{I}_{(-\infty, +\infty)}\|},
\label{errorindicator2}
\end{equation}
\added{where $x_R = x_{[\f{2N+2}{3}]}^{\b}$ is the $[\f{2N+2}{3}]^{\text{th}}$
collocation point of the generalized, $x_{0}$-shifted Hermite functions, and}
\begin{equation}
  \mathcal{F}(U_{N, x_0}^{\b}) = \frac{\|(I-\pi_{N-M, x_0}^{\b})
    U_{N, x_0}^{\b}(\cdot, t)\|}{\|U_{N, x_0}^{\b}(\cdot, t)\|},
\label{f_indicator}
\end{equation}
\added{with $M$ is taken to be $[\frac{N}{3}]$.}

\added{The major advantage of the proposed adaptive spectral method
  Fig.~\ref{algmovingscaling} is that it depends only on the numerical
  solution $U^{\beta}_{N, x_0}$ and thus does not require any prior
  knowledge on how the solution will evolve. This feature is similar
  to that of the adaptive mesh generating method which also only
  depends on the numerical solution \cite{tang2007adaptive}.  However,
  unlike the posterior error indicator that is usually used in finite
  element methods \cite{lin2014posterior}, the exterior-error and
  frequency indicators used in our adaptive spectral method does not
  directly furnish the error. The exterior-error indicator is
  specifically designed for spectral methods in unbounded domain
  problems, and controlling it by properly translating the basis
  functions can lead to a better approximation at infinity.  On the
  other hand, the frequency indicator applies to spectral methods in
  both bounded and unbounded domains, and more resembles a measure of
  the numerical error.  Ultimately, the adaptive spectral method aims
  at controlling the error by maintaining a small frequency indicator.
  While adjusting the scaling factor or changing the expansion order
  directly controls the frequency indicator, changing the displacement
  of the basis functions to control the exterior-error indicator also
  helps control the frequency indicator, as was shown in
  \cite{xia2020b}.}

Despite the numerical success of adaptive spectral methods when
applied on unbounded domains, there exists no theoretical analysis of
how the parameters $\b, x_0$, and $N$ affect the algorithm's
performance and thus far no general rule on how to best adjust these
parameters in the moving ($x_0\leftarrow \tilde{x}_0$), scaling
($\beta\leftarrow \tilde{\beta}$), and expansion order adjustment
($N\leftarrow \tilde{N}$) subroutines in order to minimize errors.
  Since the improper adjustment of $\b, x_0$, and $N$ can lead to
  large errors \citep{tang1993hermite,xiong2022short}, properly
  choosing them is crucial for the effective implementation of
  adaptive spectral methods.
\begin{table}[t!]
\vspace{-0.1in}
\scriptsize
\renewcommand*{\arraystretch}{1.25}
\begin{tabular}{| >{\centering\arraybackslash} m{9em}| 
>{\arraybackslash} m{35em}|}\hline
symbol & \hspace{3.8cm} definition \\[1pt] \hline\hline
 \,\,\, $\hat{\mathcal{H}}_{i, x_0}^{\beta}$\,\, & generalized $i^{\rm th}$-order
 Hermite function with a scaling factor $\beta$ and displacement $x_0$, 
defined in $\mathbb{R}$ as $\hat{\mathcal{H}}_{i, x_0}^{\beta}\coloneqq \hat{\mathcal{H}}_{i}(\beta(x-x_0))$ \\[1pt]  \hline
    \,\,\, $P_{N, x_0}^{\beta}$\,\, & function space 
$P_N^{x_0, \beta}\coloneqq\{\hat{\mathcal{H}}_{i, x_0}^{\beta}\}_{i=0}^N$ \\[1pt]  \hline
    \,\,\, $I$\,\, & the identity operator\\[1pt]  \hline
\,\,\, $\pi_{N, x_0}^{\beta}$\,\, & the projection operator $\pi_{N, x_0}^{\beta}$ $:L^2(\mathbb{R})\rightarrow P_N^{x_0, \beta}$ 
such that $(\pi_N^{x_0, \beta}u(x)$, $u(x) - \pi_N^{x_0, \beta}u(x))=0$ \\[1pt]  \hline
\,\,\, $\mathcal{I}_{N, x_0}^{\beta}$\,\, & the interpolation operator $\mathcal{I}_{N, x_0}^{\beta}$ $:L^2(\mathbb{R})\rightarrow P_N^{x_0, \beta}$ 
such that $\mathcal{I}_N^{x_0, \beta}u(x_i)= u(x_i)$ where $\{x_i^N\}_{i=0}^N$ are collocation points of $\{\hat{\mathcal{H}}_{i, x_0}^{\beta}\}_{i=0}^N$  \\[1pt]  \hline
   \,\,\, $U_{N, x_0}^{\beta}$\,\, & spectral expansion $U_{N, x_0}^{\beta}=\sum_{i=0}^N 
u_{i, x_0}^{\beta}\hat{\mathcal{H}}_i(\beta(x-x_0))$ \\[1pt]  \hline
  \,\,\, $ N$\,\, & expansion order of the spectral expansion \\[1pt]  \hline
  \,\,\, $ \beta$ \,\, &  scaling factor of the generalized Hermite functions \\[1pt]\hline
  \,\,\, $x_0$\,\, & displacement of the generalized Hermite functions \\[1pt]  \hline     
  \,\, $\mathcal{E}_{R}(U_{N, x_0}^{\beta}), \mathcal{E}_{L}(U_{N, x_0}^{\beta})$\,\, & $\mathcal{E}_{R}$: 
the right exterior-error indicator of the 
spectral expansion $U_{N, x_0}^{\beta}$; \qquad\quad \qquad\,  \qquad\qquad\qquad\qquad\quad $\mathcal{E}_{L}$: the left exterior-error indicator of the 
spectral expansion $U_{N, x_0}^{\beta}$ \\[1pt]  \hline
 \,\,\, $\mathcal{F}(U_{N, x_0}^{\beta})$\,\, & frequency indicator for the spectral expansion $U_{N, x_0}^{\beta}$\\[1pt]  \hline
 \,\,\, $q$\,\, & scaling factor update ($\beta$ to $\tilde{\b}$) ratio 
($\tilde{\b}\leftarrow q^n\b~\text{or}~q^{-n}\b, n\in\mathbb{N}^+$)
in the scaling technique\\[1pt]  \hline
 \,\,\, $\nu$\,\, & threshold for activating the scaling technique \\[1pt]  \hline
 \,\,\, $\delta$\,\, & minimal displacement of updating the displacement $x_0$ to $\tilde{x}_0$  ($\tilde{x}_0\leftarrow x_0+nx_0~\text{or}~x_0-nx_0, n\in\mathbb{N}^+$) in the moving technique \\[1pt]  \hline
 \,\,\, $\mu$\,\, & threshold for activating the moving technique  \\[1pt]  \hline
 \,\,\, $\eta$\,\, & threshold for increasing the number of basis functions \\[1pt]  \hline
 \,\,\, $\eta_0$\,\, & threshold for decreasing the number of basis functions \\[1pt]  \hline
  \,\,\, $\gamma$\,\, & post-refinement adjustment factor for refinement threshold 
$\tilde{\eta}\leftarrow\gamma\eta$  \\[1pt]  \hline
 \,\,\, $L^2(a, b; V)$\,\, &  space of functions $\{f:[a, b]\rightarrow{V}$ ($V$ \text{is a Banach space)} such that $f$ is measurable for $\dd{t}~ \text{and}~  \int_a^b f(t)^2\dd{t}<\infty\}$ \\[1pt]  \hline
    \,\,\, $X(t_1, t_2)$\,\, & function space $\{f: f(x, s)\in L^2((t_1, t_2), t; H^1(\mathbb{R})), \partial_s f(x, s)\in L^2((t_1, t_2),$ $ t;H^{1}(\mathbb{R}))\}$ \\[1pt]  \hline
    \,\,\, $e(t)$\,\, &  $L^2$-norm of the error
    $\|u(\cdot, t) - U_{N, x_0}^{\beta}(\cdot, t)\|_{L^2}$ at time $t$ \\[1pt]  \hline
\end{tabular}
%
\caption{\scriptsize\textbf{Overview of variables and notation.}
  List of the main variables and notations associated with the overall
  adaptive spectral method. Three key variables for adaptive spectral
  methods with generalized Hermite functions are the scaling factor
  $\b$ that determines the shape of the basis functions, the
  displacement of the basis functions $x_0$, and the expansion order
  $N$ of the spectral decomposition.}
\label{tab:model_variables}
\vspace{-0.1in}
\end{table}

In this paper, we carry out a numerical analysis of the adaptive
spectral method to specify how algorithm parameters affect the
accuracy of numerical results. We restrict ourselves to a parabolic
model problem, in any dimension, and use generalized Hermite functions
as basis functions to explore numerical performances and how
parameters in the adaptive spectral algorithm control the tuning of
the three key quantities $\b, x_0$, and $N$ in
Fig.~\ref{algmovingscaling}.  \added{Furthermore, we will explicitly
  show how the frequency indicator is related to the lower error
  bound, justifying the maintenance of a small frequency indicator in
  the adaptive spectral algorithm.}

Depending on the inverse inequality for generalized Hermite functions
\citep{Spectral2011}, such analyses for numerically solving
unbounded-domain PDEs provide a posterior error estimate. This error
estimate only relies on the numerical solution and the adjustment of
$\b, x_0$, and $N$.  Our main result is
\begin{theorem}
\rm
\label{theorem1}
The $L^2$-error at time $T$ when solving a parabolic PDE in $(x,
t)\in\mathbb{R}\times[0, T]$ with the generalized Hermite functions
and using adaptive techniques is bounded by
\begin{equation}
\begin{aligned}
e(T)\coloneqq\|u(\cdot, T) - U_{N, x_0}^{\b}(\cdot, T)\|_2\leq e_0+ e_{\text{S}} 
+ e_{\text{M}} + e_{\text{C}},
\end{aligned}
\label{errores}
\end{equation}
where $U_{N, x_0}^{\b}$ is the numerical solution; $e_0$ is the
\added{numerical discretization} error from numerically solving the
PDE. $e_{\text{S}}$ is the error bound arising from changing the
scaling factor from $\b$ to $\tilde{\b}$; $e_{\text{M}}$ is the error
bound for changing the displacement from $x_0$ to $\tilde{x}_0$;
$e_{\text{C}}$ is the error bound for coarsening, \textit{i.e.},
reducing the expansion order from $N$ to $\tilde{N}$. More
specifically, $e_{\text{S}}, e_{\text{M}}$, and $e_{\text{C}}$ take
the forms

\begin{equation}
\begin{aligned}
e_{\text{S}}&\coloneqq\sum\limits_{\rm scale} 
\frac{|\tilde{\beta}-\b|\sqrt{1+\tfrac{\tilde{\beta}}{\beta}}}{\sqrt{2}\tilde{\beta}}
\|x\partial_x U_{N, x_0}^{\beta}(\cdot, t)\|_2,\\
e_{\text{M}}&\coloneqq\sum\limits_{\rm move} 
|x_0 - \tilde{x}_0|\|\partial_x U_{N, x_0}^{\beta}(\cdot, t)\|_2,\\
e_{\text{C}}&\coloneqq \sum\limits_{\rm coarsen} 
\|(I - \pi_{\tilde{N}, x_0}^{\beta})U_{N, x_0}^{\beta}(\cdot, t)\|_2,
\end{aligned}
\label{esmc}
\end{equation}
where the sum $\sum\limits_{\rm scale}$ is taken over all scaling
steps, the sum $\sum\limits_{\rm move}$ is taken over all moving
steps, and $\sum\limits_{\rm coarsen}$ is taken over all coarsening
steps. The operators $I$ and $\pi_{\tilde{N}, x_0}^{\beta}$ are defined in 
Table~\ref{tab:model_variables}.
\end{theorem}
This result allows us to provide general guidelines for selecting the
parameters in the adaptive spectral algorithm that lead to the proper
tuning of $\b, x_0$, and $N$. \added{Specifically, the numerical
  discretization error $e_0$ in Eq.~\eqref{errores} we aim to minimize
  depends on $\b, x_0, N$. The precise dependences will be given in
  Section~\ref{error_model_pro}. Since the adaptive techniques depend
  only on the numerical solution and do not require any prior
  knowledge of the solution, the last three terms in
  Eq.~\eqref{errores} depend only on the numerical solution.  From
  this theorem, we can conclude that the smaller the adjustment in the
  scaling factor or in the displacement of the basis functions, the
  smaller the error bounds $e_{\text{S}}, e_{\text{M}}$ for carrying
  out the adaptive techniques. However, given that improper $\b$ or
  $x_0$ leads to very large $e_0$, proper dynamic adjustment of $\b$
  and $x_0$ are still needed to keep $e_0$ small, possibly at the
  expense of accumulating more error in $e_{\text{S}}, e_{\text{M}}$.}

\added{In Fig.~\ref{algmovingscaling}, the threshold
  $\mathcal{E}^m_{\text{thres}}$ is chosen to be the exterior-error
  indicator evaluated after the last adjustment of the displacement
  $x_0$, multiplied by a constant $\mu>1$.  As shown in
  \cite{xia2020b}, if the exterior-error indicator grows above such a
  threshold, the function is moving rightward, indicating that we
  should replace $x_0$ with $\tilde{x}_0 > x_0$. As
  $\lim\limits_{\tilde{x}_0\rightarrow\infty}\mathcal{E}(U_{N,
    x_0}^{\beta}(x+\tilde{x}-x_0, t))=0$, we can always find a
  $\tilde{x}_0$ such that $\mathcal{E}(U_{N,
    x_0}^{\beta}(x+\tilde{x}-x_0, t))<\mathcal{E}^m_{\text{thres}}$
  and renew $x_0\leftarrow\min\{\tilde{x}_0, x_0+d_{\max}\}$.  By the
  form of $e_{\text{M}}$ in Eq.~\eqref{esmc}, we can conclude that
  finding the smallest $\tilde{x}_0$ such that $\mathcal{E}(U_{N,
    x_0}^{\beta}(x+\tilde{x}-x_0, t))<\mathcal{E}^m_{\text{thres}}$
  while keeping $x_0-\tilde{x}_0$ small can effectively reduce
  $e_{\text{M}}$.}

\added{The scaling technique and the $p$-adaptive techniques are
  directly coupled with each other as they rely on monitoring the same
  frequency indicator.  If the function decays more slowly at
  infinity, then the frequency indicator is likely to increase,
  whereas if the function decays faster, the frequency indicator is
  likely to decrease.  When $\beta$ is to be
  decreased (more slowly decaying function), the threshold
  $\mathcal{F}^s_{\text{thres}}$ is chosen to be the frequency
  indicator after the last scaling or change of expansion order,
  multiplied by a constant $\nu>1$. When $\beta$ is to be increased (faster decaying function),
  we set the threshold to be the frequency indicator after the last
  scaling or expansion order change since a function that decreases
  more slowly is harder to approximate requiring us to be more
  tolerant of an increase in the frequency indicator. The explicit
  form of $e_{\text{S}}$ in Eq.~\eqref{esmc} suggests that to reduce
  $e_{\text{S}}$, it is desirable to find a $\tilde{\beta}$ such that
  $\tilde{\beta}-\beta$ is small.  However, there is no guarantee that
  one can find a $\tilde{\beta}$ such that $\mathcal{F}(U_{N,
    x_0}^{\tilde{\b}})<\mathcal{F}^s_{\text{thres}}$.  If the
  frequency indicator cannot be suppressed below the threshold by
  choosing $\tilde{\beta}$, a probable cause is that the function
  becomes more oscillatory, implying that the expansion order should
  be adjusted.}

\added{The $p$-adaptive threshold $\mathcal{F}^p_{\text{thres}}$ is
  chosen to be the frequency indicator after the last adjustment of
  expansion order, multiplied by a constant $\eta>1$ if refinement is
  required.  Alternatively, if coarsening is required, the threshold
  is chosen to be the frequency indicator after the last change of
  expansion order, multiplied by an another constant $\eta_0>1$ but
  $\eta_{0} < \eta$.  $\eta$ is allowed to increase with time as
  functions that oscillate rapidly are harder to approximate,
  requiring us to be more tolerant of increases in the frequency
  indicator. Since
  $\lim\limits_{\tilde{N}\rightarrow\infty}\mathcal{F}(U_{\tilde{N},
    x_0}^{\beta}(x, t))=0$, we could always find a $\tilde{N}$ such
  that $\mathcal{F}(U_{\tilde{N}, x_0}^{\beta}(x, t))\leq
  \mathcal{F}^p_{\text{thres}}$ if refinement is needed. By
  maintaining the scaling factor below the $p$-adaptive threshold
  $\mathcal{F}^p_{\text{thres}}$ and using the relationship between
  the error and the frequency indicator, the lower error bound can be
  shown to be always smaller than
  $\mathcal{F}^p_{\text{thres}}\|u(\cdot, t)\|_2 - \|(I-\pi_{N-M,
    x_0}^{\beta})u(\cdot, t)\|_2$, where $u$ is the analytical
  solution. However, tradeoffs arise.  For example, refinement itself
  does not bring about an additional error, but could result in
  additional computational cost.  On the other hand, if coarsening is
  implemented, a smaller $\tilde{N}$ could lead to a larger error
  $e_{\text{C}}$ in Eq.~\eqref{esmc} but also result in smaller
  computational cost.}



\added{In the next section, we formulate the model problem using generalized
Hermite functions and perform numerical analysis. In Section
\ref{error_adaptive}, numerical analysis for applying the adaptive
techniques is carried out and Theorem~\ref{theorem1} is proved.
Furthermore, the relationship between the error and the frequency
indicator is analyzed, explicitly explaining the efficacy of the
algorithm shown in Fig.~\ref{algmovingscaling}.  In Section
\ref{Numericalex}, numerical experiments are carried out, and an
additional improvement of the adaptive spectral method in the moving
technique is proposed and discussed. For completeness, we list the
common variables and notations in Table~\ref{tab:model_variables} that
we use throughout this paper.}

\section{Errors in solving a model problem with generalized Hermite functions}
\label{error_model_pro}
In this section, we first formulate a parabolic equation in weak form
\citep{MA2005}:

\begin{align}
\big(\p_t u(\cdot, t), v(\cdot)) + a\big(u(x, t), v(x); t\big) & 
= \big(f(\cdot, t),v(\cdot)\big),\,\,\, x\in\mathbb{R}, t\in[0, T],\,\, 
~\forall v(x)\in H^1(\mathbb{R}),
\label{model}\\
\big(u(\cdot, 0), \tilde{v}(\cdot)\big) & = \big(u_0(\cdot), 
\tilde{v}(\cdot)\big),\,\,\,  \forall v(x)\in H^1(\mathbb{R}),
\label{initial}
\end{align}
where $u_0(x)\in L^2(\mathbb{R})$ is the initial condition, $f(x, t)$
is the inhomogeneous source term (\textit{e.g.} heat source in the
heat equation), and $a(\cdot, \cdot; t)$ is a coercive symmetric
bilinear form such that there exist constants $0<c_0<C_0$ satisfying

\begin{equation}
|a(u(x, t), v(x); t)|\leq C_0\|u(\cdot, t)\|_{H^1}\, \|v(\cdot)\|_{H^1} 
\quad\text{and}\quad
c_0\|v(\cdot)\|_{H^1}^2\leq a(v, v; t), \,\,\, 
\forall u(\cdot, t), v(\cdot) \in H^1(\mathbb{R}).
\label{Acondition}
\end{equation} 
In Eqs.~\eqref{model}, \eqref{initial}, and \eqref{Acondition} and
hereafter, the inner product is taken over the spatial variable $x$,
and the norm $\|\cdot\|$ denotes the $L^2$-norm taken over $x$ unless
otherwise specified.

The solution to the model problem, Eqs.~\eqref{model} and
\eqref{initial}, exists and is unique \citep{dautray1992mathematical},
and the solution $u$ is in the so-called Bochner-Sobolev space

\begin{equation}
 W\left(0, t; H^1(\mathbb{R}), H^{-1}(\mathbb{R})\right)\coloneqq \big\{u: u(x, s)
\in L^2\left(0, t; H^1(\mathbb{R})\right), \partial_s u(x, s)
\in L^2\left(0, t;H^{-1}(\mathbb{R})\right)\big\}
\end{equation}
where $H^{-1}(\mathbb{R})$ is the dual space of $H^1(\mathbb{R})$.
For simplicity, we assume that $f(x, t)\in C(\mathbb{R}\times [0, t]),
\partial_s u(x, s)\in L^2(0, t; H^1(\mathbb{R}))$ and therefore $u\in
X(0, t)$, and its norm is given by

\begin{equation}
\|u\|_{X(0, t)}^2 = \int_0^t  \Big(\|u(\cdot, s)\|^2_{H^1(\mathbb{R})}+ 
\|\partial_s u(\cdot, s)\|^2_{H^1(\mathbb{R})}\Big)\dd{s} + \|u(\cdot, 0)\|^2.
\end{equation}

Analysis of finite element methods for solving Eqs.~\eqref{model} and
\eqref{initial} for bounded $x$ has already been performed
\citep{ueda2019inf}. Here, we wish to numerically solve
Eqs.~\eqref{model} and \eqref{initial} using spectral methods with
generalized Hermite functions.  We first fix the scaling factor $\b$, the
displacement $x_0$ of the basis functions $\hat{\mathcal{H}}_{i,
  x_0}^{\b}$, and the expansion order $N$ of the trial and test
functions. Integrating Eq.~\eqref{model} w.r.t time, we wish to
find a $U_{N, x_0}^{\beta}(x, s)\in L^2(0, t; P_{N,
  x_0}^{\beta}(\mathbb{R}))$ such that for any test function $v_{N,
  x_0}^{\b}(x, t)\in L^2(0, t; P_{N, x_0}^{\beta}(\mathbb{R}))$ and
$\tilde{v}_{N, x_0}^{\b}\in P_{N, x_0}^{\beta}(\mathbb{R})$,

\begin{equation}
\begin{aligned}
& \int_0^t \left[\big(\p_sU_{N, x_0}^{\beta}, v_{N, x_0}^{\beta}\big) 
+ a\big(U_{N, x_0}^{\beta}, v_{N, x_0}^{\beta}; t\big)\right]\dd{s} 
+\big(U_{N, x_0}^{\b}(\cdot, 0), \tilde{v}_{N, x_0}^{\beta}(\cdot)\big)\\
& \hspace{3.2cm} = \int_0^t\big(f,v_{N, x_0}^{\beta}\big)\dd{s} 
+ \big(u(\cdot, 0), \tilde{v}_{N, x_0}^{\beta}(\cdot)\big),\\
& \hspace{4.4cm} \forall (v_{N, x_0}^{\b}, \tilde{v}_{N, x_0}^{\beta}) 
\in L^2(0, t; P_{N, x_0}^{\beta}(\mathbb{R}))\times P_{N, x_0}^{\beta}(\mathbb{R}).
\label{modelpro}
\end{aligned}
\end{equation}
For notational simplicity, we denote

\begin{equation}
\textbf{v}_{N, x_0}^{\beta}\coloneqq (v_{N, x_0}^{\b}, 
\tilde{v}_{N, x_0}^{\beta}), \quad Y_{N, x_0}^{\beta}
\coloneqq L^2(0, t; P_{N, x_0}^{\beta}(\mathbb{R}))
\times P_{N, x_0}^{\beta}(\mathbb{R})\subseteq X(0, t),
\end{equation}
and equip $\textbf{v}_{N, x_0}^{\beta}\in Y_{N, x_0}^{\beta}$ with the norm

\begin{equation}
\|\textbf{v}_{N, x_0}^{\beta}\|_{Y_{N, x_0}^{\beta}}^2
=\|\big(v_{N, x_0}^{\beta}(x, t), 
\tilde{v}_{N, x_0}^{\beta}(x)\big)\|_{Y_{N, x_0}^{\beta}}^2
\coloneqq \int_0^t\|v_{N, x_0}^{\beta}(\cdot, s)\|^2_{H^1(\mathbb{R})} \dd{s} 
+ \|\tilde{v}_{N, x_0}^{\beta}(\cdot)\|^2.
\end{equation}
The solution $U_{N, x_0}^{\b}\coloneqq
\sum_{i=0}^Nu_{i,x_0}^{\b}(t)\hat{\mathcal{H}}_{i, x_0}^{\b}(x)$ of
Eq.~\eqref{modelpro} can be explicitly evaluated through the matrix
equation

\begin{equation}
\bm{u}_{N, x_0}^{\beta}(t) = e^{-{\bm A}_N^{\beta}{t}}\bm{u}_{N,
  x_0}^{\beta}(0) + e^{-{\bm A}_N^{\beta}t}\int_{0}^{t} e^{{\bm
    A}_N^{\beta}s}\bm{F}_{N, x_0}(s)\dd{s},
\label{forwardtime}
\end{equation}
where

\begin{equation}
\begin{aligned}
\bm{u}_{N, x_0}^{\beta}(s) & \coloneqq \big(u_{0, x_0}^{\b}(s),
\ldots,u_{N, x_0}^{\b}(s)\big)^T, \\
\bm{F}_{N, x_0}^{\b}(s) & \coloneqq \big(f_{0,x_0}^{\b}(s),
\ldots,f_{N, x_0}^{\b}(s)\big)^T, \\
f_{i, x_0}^{\b} & = \big(f(x, s), \hat{\mathcal{H}}^{\b}_{i, x_0}(x)\big)
\end{aligned}
\label{ufnotation}
\end{equation}
are the vectors consisting of coefficients in the spectral expansion
$U_{N, x_0}^{\b}$ and the coefficients of the spectral expansion of
the RHS term $f$ in Eq.~\eqref{modelpro}. The matrix
$\bm{A}_N^{\beta}$ is defined by

\begin{equation}
(\bm{A}_N^{\beta})_{ij}=a(\hat{\mathcal{H}}^{\b}_{i, x_0}, 
\hat{\mathcal{H}}^{\b}_{j, x_0}; t)
\label{Anotation}
\end{equation}
where $a$ is the bilinear operator in Eq.~\eqref{model}. The initial
values $u_{i, x_0}^{\b} \coloneqq \big(u(\cdot, 0),
\hat{\mathcal{H}}^{\b}_{i, x_0}(\cdot)\big)$.

Our goal is to analyze the error $e(t)=\|U_{N, x_0}^{\b}(\cdot, t)-u(\cdot,
t)\|$, where $u$ gives the solution to the model problem
(Eqs.~\eqref{model} and \eqref{initial}) and $U_{N, x_0}^{\b}$ is the
numerical solution of Eq.~\eqref{modelpro}.

\begin{theorem}
\rm
\label{theorem2}
Suppose $u$ solves Eqs.~\eqref{model} and \eqref{initial} and $U_{N,
  x_0}^{\b}$ solves Eq.~\eqref{modelpro}, then the error $e(t)=\|U_{N,
  x_0}^{\b}(\cdot, t)-u(\cdot, t)\|$ can be bounded by

\begin{equation}
e(t)\leq \frac{b_{N,\b} + B_0}{b_{N,\b}}
\inf_{z_{N, x_0}^{\beta}\in Y_{N, x_0}^{\beta}}\|u 
- z_{N, x_0}^{\beta}\|_{X(0, t)},
\label{numerror}
\end{equation}
where $B_0$ is a constant that depends on the bilinear operator
$a(\cdot, \cdot)$ and $b_{N,\b}$ depends on $a(\cdot, \cdot; t)$, the
scaling factor $\b$, and the dimension of the space $P_{N, x_0}^{\b}$.
\end{theorem}

\begin{proof}
For simplicity, we define the operator (denoting the LHS of
Eq.~\eqref{modelpro})

\begin{equation}
B(u, \textbf{v}_{N, x_0}^{\beta}) \coloneqq 
\int_0^t \Big[(\partial_s u, v_{N, x_0}^{\beta})
+ a(u, v_{N, x_0}^{\beta}; t)\Big]\dd{s} 
+ (u_0, \tilde{v}_{N, x_0}^{\beta}),\,\,\, u\in X(0, t), 
\textbf{v}_{N, x_0}^{\beta}\!\!\in Y_{N, x_0}^{\beta}.
\label{Bdef}
\end{equation}
It can be proved that $B(u, \textbf{v}_{N, x_0}^{\beta})$ is a
continuous operator, $\textit{i.e.}$, there exists a constant $B_0$
such that

\begin{equation}
B(u, \textbf{v}_{N, x_0}^{\beta})\leq 
B_0\|u\|_{X(0, t)}\|\textbf{v}_{N, x_0}^{\beta}\|_{Y_{N, x_0}^{\beta}}.
\end{equation}
Furthermore, there exists a positive constant that depends on the
dimension of the basis function space $P_{N, x_0}^{\b}$ as well as the
scaling factor $\b$ denoted by $b_{N,\b}$ such that

\begin{equation}
\inf_{0\leq U_{N, x_0}^{\beta}\in X_{N, x_0}^{\beta}}
\sup_{0\leq \textbf{v}_{N, x_0}^{\beta}\in X_{N, x_0}^{\beta}}
\frac{B(U_{N, x_0}^{\beta},\textbf{v}_{N, x_0}^{\beta})}
{\|U_{N, x_0}^{\beta}\|_{X(0, t)}\|\textbf{v}_{N, x_0}^{\beta}\|_{Y_{N, x_0}^{\beta}}}
\geq b_{N, \b}.
\label{inf_sup}
\end{equation}
Actually, we can take 
\begin{equation}
\textbf{v}_{N, x_0}^{\beta} = (U_{N,
  x_0}^{\beta}(x, s) + \frac{c_0}{(2N\b^2+1)(C_0+1)^2}\partial_sU_{N,
  x_0}^{\beta}(x, s),\, U_{N, x_0}^{\beta}(x, 0))
\label{vdef}
\end{equation}
where $c_0, C_0$ are the constants in Eq.~\eqref{Acondition}.
Therefore, by substituting $v$ as defined in Eq.~\eqref{vdef} into
Eq.~\eqref{Bdef}, we find

\begin{equation}
\begin{aligned}
B(U_{N, x_0}^{\beta}, \textbf{v}_{N, x_0}^{\beta}) \geq & 
\frac{1}{2} \big(\|U_{N, x_0}^{\beta}(\cdot, 0)\|^2 
+ \|U_{N, x_0}^{\beta}(\cdot, t)\|^2\big)\\
& \quad + c_0\int_0^t  \left(\|U_{N, x_0}^{\beta}\|^2_{H^1}  + 
\frac{1}{(2N\b^2+1)(C_0+1)^2} \|\partial_s U_{N, x_0}^{\beta}\|^2\right)\dd{s}\\
& \quad - \frac{c_0}{2}\int_0^t 
\left(\|U_{N, x_0}^{\beta}\|_{H^1}^2 
+ \frac{C_0^2}{(2N\b^2+1)^2(C_0+1)^4}\|\partial_s U_{N, x_0}^{\beta}\|_{H^1}^2\right)\dd{s}\\
\geq & \frac{1}{2} \big(\|U_{N, x_0}^{\beta}(\cdot, 0)\|^2 
+ \|U_{N, x_0}^{\beta}(\cdot, t)\|^2\big) + 
c_0\int_0^t \|U_{N, x_0}^{\beta}\|^2_{H^1}\dd{s}-
\frac{c_0}{2}\int_0^t \|U_{N, x_0}^{\beta}\|_{H^1}^2\dd{s}  \\
& \quad + \frac{c_0}{(2N\b^2+1)^2(C_0+1)^2}\int_0^t 
\|\partial_s U_{N, x_0}^{\beta}\|^2_{H^1}\dd{s} \\
& \quad -\frac{c_0}{2(2N\b^2+1)^2(C_0+1)^2}\int_0^t 
\|\partial_s U_{N, x_0}^{\beta}\|^2_{H^1}\dd{s} \\
\geq & \min\Big\{\frac{1}{2}, 
\frac{c_0}{2}, \frac{c_0}{2(2N\b^2+1)^2(C_0+1)^2}\Big\}
\|U_{N, x_0}^{\b}\|_{X(0, t)}^2 \\
\geq & \min\Big\{\frac{1}{4}, \frac{c_0}{4}, 
\frac{c_0}{2(2N\b^2+1)^2(C_0+1)^2}\Big\}\|U_{N, x_0}^{\b}\|_{X(0, t)}\;
\|\textbf{v}_{N, x_0}^{\beta}\|_{Y_{N, x_0}^{\beta}},
\end{aligned}
\end{equation}
where in the second inequality we have used the inverse inequality of
generalized Hermite functions \citep{Spectral2011} that states 

\begin{equation}
\|\partial_sU_{N,
  x_0}^{\beta}(\cdot, s)\|_{H^1}^2\leq (2N\b^2+1) \|\partial_s U_{N,
  x_0}^{\beta}(\cdot, s)\|^2.
\label{inves}
\end{equation}
Here, $b_{N, \b}\coloneqq \min\{\frac{1}{4}, \frac{c_0}{4},
\frac{c_0}{2(2N\b^2+1)^2(C_0+1)^2}\}$ is the constant that satisfies
Eq.~\eqref{inf_sup}.

For any $\textbf{v}_{N, x_0}^{\beta}\in Y_{N, x_0}^{\beta}$, if $U_{N,
  x_0}^{\beta}$ solves Eq.~\eqref{modelpro} and $u$ solves
Eqs.~\eqref{model} and \eqref{initial},

\begin{equation}
B(U_{N, x_0}^{\beta}, \textbf{v}_{N, x_0}^{\beta}) 
= B(u, \textbf{v}_{N, x_0}^{\beta}) = 
\int_0^t \big(f, v_{N, x_0}^{\beta}\big)\dd{s} 
+ (u_0, \tilde{v}_{N, x_0}^{\b}).
\label{projection}
\end{equation}
By combining Eqs.~\eqref{inf_sup} and \eqref{projection}, we find

\begin{equation}
\|U_{N, x_0}^{\beta}\|_{X(0, t)}\leq \frac{1}{b_{N, \b}}
\sup\limits_{\textbf{v}_{N, x_0}^{\beta}}\frac{B(U_{N, x_0}^{\beta}, 
\textbf{v}_{N, x_0}^{\beta})}{\|\textbf{v}_{N, x_0}^{\beta}\|_{Y_{N, x_0}^{\beta}}}
= \sup\limits_{\textbf{v}_{N, x_0}^{\beta}}
\frac{1}{b_{N, \b}}\frac{B(u, \textbf{v}_{N, x_0}^{\beta})}
{\|\textbf{v}_{N, x_0}^{\beta}\|_{Y_{N, x_0}^{\beta}}}
\leq \frac{B_0}{b_{N, \b}}\|u\|_{X(0, t)}.
\end{equation}
Finally, by the triangular inequality, we can conclude that the
approximation error is bounded:

\begin{equation}
\begin{aligned}
\|u - U_{N, x_0}^{\beta}\|_{X(0, t)}&\leq \inf_{z_{N, x_0}^{\beta}
\in Y_{N, x_0}^{\beta}}(\|u -z_{N, x_0}^{\beta}\|_{X(0, t)} 
+ \|U_{N, x_0}^{\beta} - z_{N, x_0}^{\beta}\|_{X(0, t)})\\
&\leq \frac{b_{N, \b} + B_0}{b_{N, \b}}
\inf_{z_{N, x_0}^{\beta}\in Y_{N, x_0}^{\beta}}
\|u - z_{N, x_0}^{\beta}\|_{X(0, t)}.
\end{aligned}
\end{equation}
Notice that the $L^2$-error $e(t)=\|u(\cdot, t)-U_{N, x_0}^{\beta}(\cdot,
t)\|$ at time $t$ can be bounded by $\|u - U_{N,
  x_0}^{\beta}\|_{X(0, t)}$ , and therefore Eq.~\eqref{numerror} holds. 
\end{proof}

We can also use generalized Hermite functions to numerically solve the
$D$-dimensional model problem Eq.~\eqref{modelpro},

\begin{equation}
\begin{aligned}
& \int_0^t \left[\big(\p_sU_{\bm{N}, \bm{x}_0}^{\beta}(\bm{x}, s), v_{\bm{N}, 
\bm{x}_0}^{\bm{\beta}}(\bm{x}, s)\big) + a\big(U_{\bm{N}, 
\bm{x}_0}^{\bm{\beta}}(\bm{x}, s), v_{\bm{N}, 
\bm{x}_0}^{\bm{\beta}}(\bm{x}, s; t)\big)\right]\dd{s}  \\
& \hspace{8mm}+\big(U_{\bm{N}, \bm{x}_0}^{\bm{\b}}(\bm{x}, 0), 
\tilde{v}_{\bm{N}, \bm{x}_0}^{\bm{\beta}}(\bm{x})\big) 
= \int_0^t\big(f(\bm{x}, s),v_{\bm{N}, 
\bm{x}_0}^{\bm{\beta}}(\bm{x}, s)\big)\dd{s} + 
\big(u(\bm{x}, 0), \tilde{v}_{\bm{N}, \bm{x}_0}^{\bm{\beta}}(\bm{x})\big), 
\end{aligned}
\label{modelpro2}
\end{equation}
where 

\begin{equation}
\bm{x}\coloneqq(x^1,\ldots,x^D),\,\,\,, \bm{\beta}\coloneqq(\b^1,\ldots,\b^D),\,\,\, 
\bm{x}_0\coloneqq(x_0^1,\ldots,x_0^D),\,\,\, \bm{N}\coloneqq(N^1,\ldots,N^D)
\label{multidef}
\end{equation}
are the $D$-dimensional scaling factors, displacements, and expansion
orders and

\begin{equation}
U_{\bm{N}, \bm{x}_0}^{\bm{\beta}}(\bm{x}, s), 
v_{\bm{N}, \bm{x}_0}^{\bm{\beta}}(\bm{x}, s)
\in L^2\big(0, t; \bigotimes_{h=1}^D P_{N^h, x_0^h}^{\beta^h}(\mathbb{R})\big), 
\quad \tilde{v}_{\bm{N}, \bm{x}_0}^{\bm{\beta}}\in 
\bigotimes_{h=1}^D P_{N^h, x_0^h}^{\beta^h}(\mathbb{R}).
\end{equation} 
A multiple dimension version of the error bound Eq.~\eqref{numerror}
can be similarly derived

\begin{equation}
\begin{aligned}
\|u(\cdot, t) - U_{\bm{N}, \bm{x}_0}^{\bm{\beta}}(\cdot, t)
\|\leq \frac{b_{\bm N, \bm \b} + B_0}{b_{\bm N, \bm \b}}
\inf_{z_{\bm{N}, \bm{x}_0}^{\bm \beta}\in {Y}_{\bm{N}, 
\bm{x}_0}^{\bm{\beta}}}\|u - z_{{\bm N}, \bm{x}_0}^{\bm{\beta}}\|_{X(0, t)},
\end{aligned}
\label{erroresti_multi}
\end{equation}
where $b_{\bm N, \bm \b}\coloneqq \min\{\frac{1}{4}, \frac{c_0}{4},
\frac{c_0}{2(\sum_{h=1}^D2N_i\b_i^2+1)(C_0+1)^2}\}$. The function
spaces are
\begin{equation}
\begin{aligned}
X(0, t)&\coloneqq \Big\{u: u(\bm{x}, s)\in L^2\left(0, t; H^1(\mathbb{R}^D)\right), 
\partial_s u(\bm{x}, s)\in L^2\left(0, t;H^{1}(\mathbb{R}^D)\right)\Big\}.\\
{Y}_{\bm{N}, \bm{x}_0}^{\bm{\beta}}&\coloneqq L^2\big(0, t; 
\bigotimes_{h=1}^D P_{N^h, x_0^h}^{\beta^h}(\mathbb{R})\big)
\times\bigotimes_{h=1}^D P_{N^h, x_0^h}^{\beta^h}(\mathbb{R}).
\end{aligned}
\end{equation}

\section{Errors of adaptive techniques}
\label{error_adaptive}
In this section, we analyze the errors directly associated with the
moving, scaling, and $p$-adaptive techniques that automatically change
the shape, the translation, and the order of the numerical solution
through adjustment of $\beta$, $x_0$, and $N$, respectively
\citep{xia2020b,xia2020a}.  We derive the error bound when solving
Eq.~\eqref{modelpro} and prove Theorem~\ref{theorem1} presented in
Introduction.  Doing so explicitly shows how changing $\beta, x_0$,
and $N$ affects the error, thus providing insight on how to choose
parameters in the adaptive algorithm that leads to the proper tuning of
$\beta, x_0$, and $N$.

Instead of using collocation methods to carry out the scaling, moving,
or $p$-adaptive methods as was done in previous work
\citep{xia2020b,xia2020a} (\textit{i.e.}, enforcing the updated
numerical solution to be the same with the original numerical solution
on the new collocation points), we now use the Galerkin method
(\textit{i.e.}, projecting the numerical
solution onto the space of adjusted basis functions).  For example,
given the numerical solution $U_{N, x_0}^{\beta}(x, t)$ at time $t$,
if we change its scaling factor from $\b$ to $\tilde{\b}$, previous
implementation in \citep{xia2020b,xia2020a} replaces $U_{N,
  x_0}^{\beta}(x, t)$ with $\mathcal{I}_{N, x_0}^{\tilde{\b}}U_{N,
  x_0}^{\beta}\in P_{N, x_0}^{\tilde{\b}}$ as the new numerical
solution. This new numerical solution $\mathcal{I}_{N,
  x_0}^{\tilde{\b}}U_{N, x_0}^{\beta}$ takes on the same values as
$U_{N, x_0}^{\beta}$ at the collocation points for the new basis
functions $\{\hat{\mathcal{H}}_{i,
  x_0}^{\tilde{\b}}\}_{i=0}^N$. Therefore, the error after changing
$\b$ to $\tilde{\b}$ and replacing $U_{N, x_0}^{\beta}(x, t)$ with
$\mathcal{I}_{N, x_0}^{\tilde{\b}}U_{N, x_0}^{\beta}$ can be bounded
by
\begin{equation}
\|u - \mathcal{I}_{N, x_0}^{\tilde{\b}}U_{N, x_0}^{\beta}\| \leq 
\|u - U_{N, x_0}^{\beta}\|+
\|(I-\mathcal{I}_{N, x_0}^{\tilde{\b}})U_{N, x_0}^{\beta}\|.
\label{interp}
\end{equation}
In this work, we project the numerical solution onto $P_{N,
  x_0}^{\tilde{\b}}\coloneqq\{\hat{\mathcal{H}}_{i,
  x_0}^{\tilde{\b}}\}_{i=0}^N$, \textit{i.e.}, using $\pi_{N,
  x_0}^{\tilde{\b}}U_{N, x_0}^{\beta}$ as the new numerical
solution. Therefore, the error bound after changing the scaling factor
is
\begin{equation}
\|u - \pi_{N, x_0}^{\tilde{\b}}U_{N, x_0}^{\beta}\|\leq 
\|u - U_{N, x_0}^{\beta}\|+\|(I-\pi_{N, x_0}^{\tilde{\b}})U_{N, x_0}^{\beta}\|.
\label{projp}
\end{equation}
The second term on the RHSs of Eqs.~\eqref{interp} and~\eqref{projp}
can be viewed as an additional error bound resulting from changing the
scaling factor. Furthermore, we are able to show
\begin{equation}
\|(I-\mathcal{I}_{N, x_0}^{\tilde{\b}})U_{N,
  x_0}^{\beta}\|\geq \|(I-\pi_{N,
  x_0}^{\tilde{\b}})U_{N, x_0}^{\beta}\|.
\label{proj_prop}
\end{equation}
The proof is straightforward.  Assuming the spectral expansion of
$U_{N, x_0}^{\beta}$ under the new basis functions
$\{\hat{\mathcal{H}}_{i, x_0}^{\beta}\}$ is

\begin{equation}
U_{N, x_0}^{\beta}=\sum_{i=0}^{\infty}u_{i, x_0}^{\beta}\hat{\mathcal{H}}_{i, x_0}^{\beta}.
\end{equation}
By definition, 
\begin{equation}
\begin{aligned}
\pi_{N, x_0}^{\tilde{\b}}U_{N, x_0}^{\beta} = \sum_{i=0}^{N}u_{i, x_0}^{\beta}
\hat{\mathcal{H}}_{i, x_0}^{\beta},
\quad \mathcal{I}_{N, x_0}^{\tilde{\b}}U_{N,x_0}^{\beta} 
=  \sum_{i=0}^{N}\tilde{u}_{i, x_0}^{\beta}\hat{\mathcal{H}}_{i, x_0}^{\beta}.
\end{aligned}
\end{equation}
Therefore, 

\begin{equation}
\begin{aligned}
\|(I-\mathcal{I}_{N, x_0}^{\tilde{\b}})U_{N, x_0}^{\beta}\| 
& = \left[\sum_{i=0}^{N}(\tilde{u}_{i, x_0}^{\beta}-u_{i, x_0}^{\beta})^2
\|\hat{\mathcal{H}}_{i, x_0}^{\beta}\|^2 + \sum_{i=N+1}^{\infty}(u_{i, x_0}^{\beta})^2
\|\hat{\mathcal{H}}_{i, x_0}^{\beta}\|^2 \right]^{\tfrac{1}{2}}\\
\: &\geq \left[\sum_{i=N+1}^{\infty}(u_{i, x_0}^{\beta})^2
\|\hat{\mathcal{H}}_{i, x_0}^{\beta}\|^2\right]^{\tfrac{1}{2}} = 
\|(I-\pi_{N,x_0}^{\tilde{\b}})U_{N, x_0}^{\beta}\|.
\end{aligned}
\end{equation}
With Eq.~\eqref{proj_prop}, using the projected $\pi_{N,
  x_0}^{\tilde{\b}}U_{N, x_0}^{\beta}$ as the new numerical solution
instead of the interpolated $\mathcal{I}_{N, x_0}^{\tilde{\b}}U_{N,
  x_0}^{\beta}$ might lead to a smaller error bound.

\subsection{Posterior error estimate}

We derive the posterior error estimates that depend on the numerical
solution $U_{N, x_0}^{\b}\in P_{N, x_0}^{\b}$ and on how $\b, x_0$,
and $N$ are changed. Combining the error estimate of the adaptive
techniques with Theorem~\ref{theorem2}, the error estimate for
numerically solving Eqs.~\eqref{model} and \eqref{initial}, our
ultimate goal is to prove Theorem~\ref{theorem1}, the error estimate
for adaptive spectral methods. To start, we analyze the errors from
the three adaptive techniques.

\subsubsection{Scaling technique error}

\added{First, we derive the error bound associated with changing the
  scaling factor $\b$, which corresponds to the scaling technique
  error $e_{\text{S}}$ in Eq.~\eqref{errores} of
  Theorem~\ref{theorem1}.} Suppose at time $t$, we change $\b$ to
  $\tilde{\b}$ and replace the numerical solution $U_{N, x_0}^{\b}$
  with $\pi_{N, x_0}^{\tilde{\beta}}U_{N, x_0}^{\b}\in P_{N,
    x_0}^{\tilde{\b}}$, the error is

\begin{equation}
\begin{aligned}
\|u(\cdot, t)-\pi_{N, x_0}^{\tilde{\beta}}U_{N, x_0}^{\beta}(\cdot, t)\| \leq 
\|u(\cdot, t) - U_{N, x_0}^{\beta}(\cdot, t)\| + \|(
I-\pi_{N, x_0}^{\tilde{\beta}})U_{N, x_0}^{\beta}(\cdot, t)\|
\end{aligned}
\end{equation}
where the first term on the RHS is the error before scaling and the
second term on the RHS is the additional error bound from changing the
scaling factor (``scaling error''). Denoting $\b' = \tilde{\b}/\b$, we
can further bound the scaling error by

\begin{equation}
\begin{aligned}
\|(I-\pi_{N, x_0}^{\tilde{\beta}})U_{N, x_0}^{\beta}(\cdot, t)\| & \leq 
\|U_{N, x_0}^{\beta}(x, t) - U_{N, x_0}^{\beta}(\b' x, t)\| \\
 & =\left[\int_{\mathbb{R}} \left(\int_{\b'{x}}^x 
\partial_y{U}_{N, x_0}^{\beta}(y, t)\dd{y}\right)^2\dd{x}\right]^{\frac{1}{2}} \\
 & \leq \left[\int_{\mathbb{R}}|1-\b'|x \left(\int_{\b'{x}}^x 
\Big(\partial_y{U}_{N, x_0}^{\beta}(y, t)\Big)^2\dd{y}\right)\dd{x}\right]^{\frac{1}{2}} \\
 & = \frac{|1-\b'|\sqrt{1+\beta'}}{\sqrt{2}\beta'}
\|x\partial_x{U}_{N, x_0}^{\beta}(x, t)\|.
\end{aligned}
\end{equation}
Therefore, the error after changing the scaling factor from $\b$ to
$\tilde{\b}$ is bounded by

\begin{equation}
\|u(\cdot, t)-\pi_{N, x_0}^{\tilde{\beta}}U_{N, x_0}^{\beta}(\cdot, t)\|
\leq \|u(\cdot, t) - U_{N, x_0}^{\beta}(\cdot, t)\| + 
\frac{|1-\b'|\sqrt{1+\beta'}}{\sqrt{2}\beta'}\|x\partial_x{U}_{N, x_0}^{\beta}(x, t)\|.
\label{scaleerror}
\end{equation}
From Eq.~\eqref{scaleerror}, the second term in the last equality is
the additional error bound resulting from scaling. The factor
$\frac{|1-\b'|\sqrt{1+\beta'}}{\sqrt{2}\beta'}$ is directly related to
how much the scaling factor is changed while $\|x\partial_x{U}_{N,
  x_0}^{\beta}(x, t)\|$ depends on the spatial derivative of the
pre-scaled solution.

\subsubsection{Moving technique error}

\added{Next, we derive the error bound associated with changing the
  displacement $x_0$, which corresponds to the moving technique error
  $e_{\text{M}}$ in Eq.~\eqref{errores} of Theorem~\ref{theorem1}.}
Given the numerical solution $U_{N, x_0}^{\b}$, if we change the
displacement of the basis functions from $x_0$ to $\tilde{x}_0$ and
set $\pi_{N, \tilde{x}_0}^{\b}U_{N, \tilde{x}_0}^{\b}\in P_{N,
  \tilde{x}_0}^{\b}$ as the new numerical solution, the error is

\begin{equation}
\begin{aligned}
\|u(\cdot, t)-\pi_{N, \tilde{x}_0}^{\b}U_{N, \tilde{x}_0}^{\b}(\cdot, t)\|
\leq \|u(\cdot, t) - U_{N, x_0}^{\b}(\cdot, t)\| 
+ \|(I-\pi_{N, \tilde{x}_0}^{\b})U_{N, x_0}^{\b}(\cdot, t)\|,
\end{aligned}
\end{equation}
where the second term on the RHS is the additional error bound from
changing $x_0$ (``moving error''). Furthermore, it is bounded by

\begin{equation}
\begin{aligned}
\|(\pi_{N, \tilde{x}_0}^{\b}-I)U_{N, x_0}^{\b}(\cdot, t)\| &
\leq \|U_{N, x_0}^{\beta}(x, t) - U_{N, x_0}^{\beta}(x-\tilde{x}_0+x_0, t)\| \\
& \leq  \left[\int_{\mathbb{R}}|x_0-\tilde{x}_0|
\left(\int_{x-\tilde{x}_0+x_0}^{x}\big(\partial_y{U}_{N, x_0}^{\beta}(y, t)\big)^2 
\dd{y}\right) \dd{x}\right]^{\tfrac{1}{2}} \\
& = d\| \partial_x{U}_{N, x_0}^{\beta}(\cdot, t)\|,
\end{aligned}
\label{moving}
\end{equation}
where $d \coloneqq |x_0-\tilde{x}_0|$. Thus, the error after
changing the displacement from $x_0$ to $\tilde{x}_0$ is bounded by

\begin{equation}
\|u(\cdot, t)-\pi_{N, x_0}^{\beta}U_{N, \tilde{x}_0}^{\beta}(\cdot, t)\|
\leq \|u(\cdot, t) - U_{N, x_0}^{\b}(\cdot, t)\|
+d\| \partial_x{U}_{N, x_0}^{\beta}(\cdot, t)\|.
\end{equation}
We see that the additional error bound associated with moving depends
on the change in the displacement $x_0$ and the spatial derivative
$\partial_x{U}_{N, x_0}^{\beta}(x, t)$ of the pre-translated numerical
solution.

\subsubsection{$p$-adaptive technique error}
Finally, we analyze the error associated with the $p$-adaptive
technique, \added{which corresponds to the $p$-adaptive technique
  error $e_{\text{C}}$ in Eq.~\eqref{errores} of
  Theorem~\ref{theorem1}}. When projecting the numerical solution
$U_{N, x_0}^{\b}$ onto the new space $P_{\tilde{N},x_0}^{\b}$, no
extra error will be introduced when $\tilde{N}>N$ (refinement) because
the basis functions $\{\hat{\mathcal{H}}_{i,
  x_0}^{\b}\}_{i=0}^{\tilde{N}}$ form an orthogonal set of basis
functions and $\pi_{\tilde{N}, x_0}^{\b}U_{N, x_0}^{\b} = U_{N,
  x_0}^{\b}$, \textit{i.e.},

\begin{equation}
\|u(\cdot, t)-\pi_{\tilde{N}, x_0}^{\b}U_{N, x_0}^{\b}(\cdot, t)\|
=\|u(\cdot, t) - U_{N, x_0}^{\b}(\cdot, t)\|,\, ~\tilde{N}>N.
\end{equation}
When we reduce the number of basis functions from $N$ to $\tilde{N}<N$
(coarsening), we use $\pi_{\tilde{N}, x_0}^{\b}U_{N, x_0}^{\b} =
U_{\tilde{N}, x_0}^{\b}$ as the new numerical
solution. $\pi_{\tilde{N}, x_0}^{\b}U_{N, x_0}^{\b}$ leaves out the
last $N-\tilde{N}$ terms in the spectral expansion of
$U_{N,x_0}^{\b}$. Therefore, the error after coarsening can be bounded
by

\begin{equation}
\begin{aligned}
\|u(\cdot, t) - \pi_{\tilde{N}, x_0}^{\b}U_{N, x_0}^{\b}(\cdot, t)\|
\leq\|u(\cdot, t)-U_{N, x_0}^{\b}(\cdot, t)\|+
\|(I-\pi_{\tilde{N}, x_0}^{\b})U_{N, x_0}^{\b}(\cdot, t)\|,\, ~\tilde{N}<N.
\end{aligned}
\label{refinement}
\end{equation}
In Eq.~\eqref{refinement}, the second term in the last inequality is
the additional error bound that results from truncating the spectral
expansion and leaving out the last $N-\tilde{N}$ terms.

Next, we generalize Theorem~\ref{theorem2} to forward time from
$t_{0}$ to $t_{1}$ given $U_{N, x_0}^{\b}(x, t_0)$. We assume that no
adaptive technique is activated within $t\in(t_0, t_{1})$ and denote
$e(x, t) = u(x, t)-U_{N, x_0}^{\b}(x, t), t\in[t_0, t_1]$, where $u$
is the solution to Eqs.~\eqref{model} and \eqref{initial}.  The error
at $t_1$, $e(x, t_{1}) = u(x, t_1)-U_{N, x_0}^{\b}(x, t_1)$, can be
decomposed as $e(x, t_{1}) = e_1(x, t_{1}) + e_2(x, t_{1})$ where
$e_1(x, t_{1})$ is the error $u(x,t_1) -
\tilde{U}_{N,x_0}^{\b}(x,t_1)$ with $\tilde{U}_{N, x_0}^{\b}$ solving
Eq.~\eqref{modelpro} with initial condition $u(x, t_0)$. The second
error term $e_2(x, t_1)\in L^2(t_0, t_1; P_{N, x_0}^{\b})$ satisfies
\begin{equation}
\begin{aligned}
& \int_{t_0}^{t_{1}} \left(\p_s e_2, v\right) 
 + a\left(e_2, v; t\right)\dd{s} + \left(e_2(\cdot, t_0), \tilde{v}(\cdot, t_0)\right)\\
& \qquad \quad=  \left(e(\cdot, t_0), \tilde{v}(\cdot, t_0)\right), \quad  \forall\,  
v\in L^2(t_0, t_1;  
P_{N, x_0}^{\b}),\,\, \tilde{v}\in P_{N, x_0}^{\b}.
\end{aligned}
\label{e2eqn}
\end{equation}
From Theorem~\ref{theorem2},

\begin{equation}
\|e_1(\cdot, t_1)\| \leq \frac{b_{N, \b}+B_0}{b_{N, \b}}
\, \|(I-\pi_{N, {x_0}}^{\beta})u\|_{X(t_0, t_{1})}.
\end{equation}
Additionally, since the bilinear form $a(\cdot, \cdot)$ is positive
definite, substituting $v(x, t)=e_2(x, t)$ and $\tilde{v}=e_2(x, t_i)$
into Eq.~\eqref{e2eqn}, we conclude that $\|e_2(\cdot, t_{1})\| \leq
\|e(\cdot,t_{0})\|=e(t_0)$.  Therefore,

\begin{equation}
\begin{aligned}
e(t_{1}) \leq e( t_0) + \frac{b_{N, \b}+B_0}{b_{N, \b}}
\, \|(I-\pi_{N, {x_0}}^{\beta})u\|_{X(t_i, t_{i+1})}.
\end{aligned}
\label{errorest}
\end{equation}
Specifically, this error bound does not depend on the step size
$\Delta{t} = t_{i+1}-t_i$ if we use

\begin{equation}
\bm{u}_{N, x_0}^{\beta}(t+\Delta{t}) = e^{-{\bm A}_N^{\beta}\Delta{t}}
\bm{u}_{N, x_0}^{\beta}(t) +  e^{-{\bm A}_N^{\beta}\Delta{t}}
\int_{t}^{t+\Delta{t}} e^{{\bm A}_N^{\beta}(s-t)}\bm{F}_{N, x_0}(s)\dd{s},
\label{time_forward}
\end{equation}
with $\bm{u}_{N, x_0}^{\beta},\, \bm{F}_{N, x_0}^{\beta}$ defined by
Eq.~\eqref{ufnotation} and $\bm{A}_N^{\beta}$ defined by
Eq.~\eqref{Anotation}.

Now, we are ready to prove Theorem~\ref{theorem1}, the overall error
bound using the adaptive spectral methods.  We define the times of the
$\ell^{\text{th}}$ scaling, the $\ell^{\text{th}}$ \added{moving}, and
the $\ell^{\text{th}}$ changing \added{of} the expansion order to be
$t^s_{\ell}, t^m_{\ell}$, and $t_{\ell}^c$, respectively.  We denote
the scaling factors right before the $\ell^{\text{th}}$ scaling,
moving, and changing the expansion order to be $\b_{\ell}^s,
\b_{\ell}^m$, and $\b_{\ell}^c$, the displacements right before the
$\ell^{\text{th}}$ scaling, moving, and changing the expansion order
to be ${x_0}_{\ell}^s, {x_0}_{\ell}^m$, and ${x_0}_{\ell}^c$, and the
expansion orders right before the $\ell^{\text{th}}$ scaling, moving,
and changing the expansion order to be $N_{\ell}^s, N_{\ell}^m$, and
$N_{\ell}^c$, respectively.  After the $\ell^{\text{th}}$ scaling, we
denote the new scaling factor to be $\tilde{\b}_{\ell}^s$ and the
ratio ${\b'}^s_{\ell}\coloneqq \tilde{\b}_{\ell}^{s}/\b_{\ell}^s$;
after the $\ell^{\text{th}}$ moving, we denote the new displacement to
be ${\tilde{x_0}}_{\ell}^m$ and $d_{\ell}^m\coloneqq
|{\tilde{x_0}}_{\ell}^m - {x_0}_{\ell}^m|$; after the
$\ell^{\text{th}}$ change of the expansion order, we denote the new
expansion order as $\tilde{N}_{\ell}^c$.

The times at which the scaling factor or the displacement of the basis
functions is changed, or the expansion order is reduced, are indicated
by $t_i$ in chronological order $0 = t_0 \leq t_1...\leq t_i\leq
t_{K^s+K^m+K^c+1}=T$, where $K^s$, $K^m$, and $K^c$ are the total
number of scalings, movings, and changing the expansion order
within $t\in[0, T]$. Specifically, if $t_i=t_{i+1}$, then more than
one adaptation is triggered simultaneously. The corresponding constant
that satisfies the inequality Eq.~\eqref{numerror} during $[t_i,
  t_{i+1}]$ is denoted as $(b_{N_i, \b_i}+ B_0)/b_{N_i,\b_i}$. From
the error estimates of the scaling, moving, and $p$-adaptive
techniques in Eqs.~\eqref{scaleerror}, \eqref{moving},
\eqref{refinement}, and Eq.~\eqref{errorest}, we conclude

\begin{equation}
\begin{aligned}
e(T)\leq & \sum_{i=0}^{K^s+K^m+K^c}\frac{b_{N_i, \b_i}+B_0}{b_{N_i, \b_i}}
\,\|(I-\pi_{N_i, {x_0}_i}^{\beta_i})u\|_{X(t_i, t_{i+1})}\\
& \qquad + \sum_{\ell=1}^{K^s}\frac{|1-{\b'}^s_{\ell}|
  \sqrt{1+{\beta'}^s_{\ell}}}{\sqrt{2}{\beta'}^s_{\ell}}\,
\|x{\partial_xU}_{N_{\ell}^s, {x_{0}}_{\ell}^s}^{\beta_{\ell}^s}(x, t_{\ell}^s)\|\\
& \qquad + \sum_{\ell=1}^{K^m}d_{\ell}^m \|\partial_xU_{N_{\ell}^m, 
{x_0}_{\ell}^m}^{\beta_{\ell}^m}(\cdot, t_{\ell}^m)\|\\
& \qquad +\sum_{\ell=1}^{K^c}\|(I-\pi_{\tilde{N}_{\ell}^c, {x_0}_{\ell}^c}^{\b_{\ell}^c})
U_{N_{\ell}^c, {x_0}_{\ell}^c}^{\b_{\ell}^c}(\cdot, t_{\ell}^c)\|\\
\leq & \sum_{i=0}^{K^s+K^m+K^c} 
\frac{b_{N_i, \b_i}+B_0}{b_{N_i, \b_i}}\,
\|(I-\pi_{N_i, {x_0}_i}^{\beta_i})u\|_{X(t_i, t_{i+1})}\\
& \qquad + \sum_{\ell=1}^{K^s}\frac{|1-{\b'}^s_{\ell}|
\sqrt{1+{\beta'}^s_{\ell}}}{\sqrt{2}{\beta'}^s_{\ell}}
(2N_{\ell}^s+1)\,\|{U}_{N_{\ell}^s, {x_{0}}_{\ell}^s}^{\beta_{\ell}^s}(\cdot, t_{\ell}^s)\|\\
& \qquad + \sum_{\ell=1}^{K^m}\sqrt{(2N_{\ell}^m+1)}\beta_{\ell}^md_{\ell}^m\,
\|U_{N_{\ell}^m,{x_0}_{\ell}^m}^{\beta_{\ell}^m}(\cdot, t_{\ell}^m)\|\\
& \qquad + \sum_{{\ell}=1}^{K^c}
\|(I-\pi_{\tilde{N}_{\ell}^c,{x_0}_{\ell}^c}^{\b_{\ell}^c})
U_{N_{\ell}^c, {x_0}_{\ell}^c}^{\b_{\ell}^c}(\cdot, t_{\ell}^c)\|
\end{aligned}
\label{errorbound}
\end{equation}
where we have used the three-term recurrence relation for
generalized Hermite functions and the inverse inequality
Eq.~\eqref{inves} to bound $\|x{\partial_xU}_{N_{\ell}^s,
  {x_{0}}_{\ell}^s}^{\beta_{\ell}^s}(x, t_{\ell}^s)\|$ and
$\|\partial_xU_{N_{\ell}^m, {x_0}_{\ell}^m}^{\beta_{\ell}^m}(\cdot,
t_{\ell}^m)\|$ in the second inequality. Note that in the first term
of Eq.~\eqref{errorbound}, if $t_i=t_{i+1}$ then we define
$\|(I-\pi_{N_i, {x_0}_i}^{\beta_i})u\|_{X(t_i, t_{i+1})}\coloneqq 0$.
The first term on the RHS of last inequality corresponds to $e_0$ in
Theorem~\ref{theorem1}, and the second, third, and last terms on the
RHS of last inequality correspond to $e_{\rm S}, e_{\rm M}$, and
$e_{\rm C}$, respectively. \added{Note that in Eq.~\eqref{errorbound},
  the first, second, third, and fourth terms on the RHS give the exact
  forms of $e_0, e_{\text{S}}, e_{\text{M}}$, and $e_{\text{C}}$ in
  Eq.~\eqref{errores} of Theorem~\ref{theorem1}.}

From Eq.~\eqref{errorbound}, the errors caused by scaling and moving
(the second and third terms of the equation) suggest that the smaller
the adjustment in $\b$ or $x_0$, the smaller the factors
$|1-{\b'}_{\ell}^s|$ and $d_{\ell}^m$ in the scaling or moving errors.
Therefore, we should set the triggering parameters $q \lesssim 1$
($\lesssim$ means smaller but close to) and $0 \lesssim \delta$ in
Table~\ref{tab:model_variables} so that the scaling factor $\b$ and
the displacement $x_0$ can be tuned more accurately without
over-adjustment that may lead to larger errors.

When coarsening, decreasing the expansion order $N$ too much will
increase the coarsening error through the last term in
Eq.~\eqref{errorbound}. Increasing the coarsening threshold $\eta_0$
to make it harder to decrease $N$ can preserve accuracy but possibly
at the expense of keeping a higher computational burden. Note that
although the effect of refinement does not explicitly reveal itself in
the error bound Eq.~\eqref{errorbound}, both a smaller initial
refinement threshold $\eta$ and a smaller $\gamma$ (the ratio of
increasing the refinement threshold) could lead to larger $N$ and thus
smaller errors (the first term of the second equation in
Eq.~\eqref{errorbound}).  However, if $N$ increases, so will the
computational cost. Using the numerical example presented in the next
section, we will discuss how to set $\gamma$ and $\eta$ so that high
accuracy can be achieved without significant degradation of
computational efficiency.  Since the adaptive techniques do not
require prior information on the solution, the last three terms in
Eq.~\eqref{errorbound}, \textit{i.e.}, errors from adaptive
techniques, depend only on the latest numerical solution itself.

Note that the numerical error in solving Eqs.~\eqref{model} and
\eqref{initial} is no less than the projection error

\begin{equation}
e(T)=\|u(\cdot, t) - U_{N, x_0}^{\beta}(\cdot, t)\|\geq 
\|(I-\pi_{N, x_0}^{\beta})u(\cdot, t)\|,
\label{lowerbound}
\end{equation}
and it has also been shown that improper scaling of generalized
Hermite functions can lead to large projection errors
\citep{tang1993hermite}. Furthermore, in Examples 2, 3, 5 in
\citep{xia2020a} and Example 2 in \citep{xia2020b}, improper
displacement $x_0$ or a too-small expansion order $N$ will also lead
to projection errors, implying a large $e(T)$. Therefore, timely and
accurate implementation of the adaptive techniques is important for
controlling the lower error bound (the projection error)
Eq.~\eqref{lowerbound}.  Consequently, to adjust them properly, we
need to set $1 \lesssim \nu$ and $1 \lesssim \mu$ in the scaling and
moving technique algorithms, respectively.

A $D$-dimensional generalization of Eq.~\eqref{errorbound} for spatial
variables $\bm{x}=(x_1,...,x_D)\in\mathbb{R}^D$ can be similarly
derived using $U_{\bm{N}, {\bm{x_0}}}^{\bm{\beta}}(\bm{x}, t)\coloneqq
\sum_{i_1=0}^{N^1}...\sum_{i_D=0}^{N^D} u_{i^1,...,i^D,
  \bm{x}_0}^{\bm{\b}}(t)\Pi_{h=1}^D\hat{\mathcal{H}}_{i^h,x_{0}^{h}}^{\b^h}(x)$:
\begin{equation}
\begin{aligned}
e(T)= & \|u(\cdot, t) -U_{\bm{N}, {\bm{x}_0}}^{\bm{\beta}}(\cdot, t)\| \\
\leq & \sum_{i=0}^{\bm{K}^s+\bm{K}^m+\bm{K}^c} 
\frac{b_{\bm{N}_i, \bm{\b}_i}+B_0}{b_{\bm{N}_i, \bm{\b}_i}}
\|(I-\pi_{\bm{N}_i, {\bm{x}_0}_i}^{\bm{\beta}_i})u\|_{X(t_i, t_{i+1})}\\
& \qquad + \sum_{h=1}^D\sum_{\ell=1}^{K^{h, s}}\frac{|1-{\b'}^{h, s}_{\ell}|
\sqrt{1+{\beta'}^{h, s}_{\ell}}}{\sqrt{2}{\beta'}^{h, s}_{\ell}}
(2N_{\ell}^{h, s}+1)
\|U_{\bm{N}_{\ell}^{h, s}, {{\bm{x}_0}^{h, s}_{{\ell}}}}^{\bm{\beta}_{\ell}^{h, s}}
(\cdot, t_{\ell}^{h, s})\| \\
& \qquad + \sum_{h=1}^D\sum_{\ell=1}^{K^{h, m}}
\sqrt{2N_{\ell}^{h, m}+1}{\beta}_{\ell}^{h, m}d_{\ell}^{h, m}
\|U_{\bm{N}_{\ell}^{h, m}, 
{\bm{x}_0}_{\ell}^{h, m}}^{\bm{\beta}_{\ell}^{h, m}}(\cdot, t_{\ell}^{h, m})\|\\
& \qquad + \sum_{h=1}^D\sum_{\ell=1}^{K^{h, c}}
\|(I - \pi_{{{\bm{\tilde{N}}}_r}^{h, c}, {{\bm{x}_0}^{h, c}_{\ell}}}^{\bm{\b}_r^{h, c}})
U_{{\bm{N}_r}^{h, c}, {{\bm{x}_0}^{h, c}_{\ell}}}^{\bm{\b}_r^{h, c}}(\cdot, t_{\ell}^{h, c})\|
\end{aligned}
\label{multidimension}
\end{equation}
where $\bm{\beta}, {\bm{x_0}}$, and $\bm{N}$ are the corresponding
$D$-dimensional scaling factor, displacement, and expansion order
defined in Eq.~\eqref{multidef}. $\bm{K}^s = \sum_{h=1}^D K^{h, s},
\bm{K}^m = \sum_{h=1}^DK^{h, m}, \bm{K}^c = \sum_{h=1}^DK^{h, c}$ are
the total number of times of performing scaling, moving, and changing
the expansion orders, across all dimensions ($K^{h, s}, K^{h, m},
K^{h, c}$ are the numbers of using the scaling, moving, or
$p$-adaptive technique in the $h^{\text{th}}$ dimension,
respectively), the constant $(b_{\bm{N}_i,
  \bm{\b}_i}+B_0)/b_{\bm{N}_i, \bm{\b}_i}$ is the RHS constant in the
inequality~\eqref{erroresti_multi} during $[t_i, t_{i+1}]$, and
$t_{\ell}^{h, s}, t_{\ell}^{h, m}, t_{\ell}^{h, c}$ are the times of
the $\ell^{\text{th}}$ scaling, moving, or changing the expansion
order in the $h^{\text{th}}$ dimension, respectively. The second,
third and last terms in Eq.~\eqref{multidimension} describe scaling
error bounds in all dimensions, moving error bounds in all dimensions,
and coarsening error bounds in all dimensions.

In Eq.~\eqref{multidimension}, $\bm{\beta}^{h, s}_{\ell}\coloneqq
({\beta}^{1, s}_{\ell},...,{\beta}^{D, s}_{\ell}), \bm{\beta}^{h,
  m}_{\ell}$, and $\bm{\beta}^{h, c}_{\ell}$ are the $D$-dimensional
scaling factors right before the $\ell^{\text{th}}$ scaling, moving,
or changing the expansion order in the $h^{\text{th}}$
dimension. Similarly, ${\bm{x}_0}_{{\ell}}^{h,
  s}\coloneqq({{x}_0}_{\ell}^{1, s},..., {x_0}_{\ell}^{D, s}),
{\bm{x}_0}_{\ell}^{h, m}$, and ${\bm{x}_0}_{\ell}^{h, c}$ are the
$D$-dimensional displacements right before the $\ell^{\text{th}}$
scaling, moving, or change of expansion order in the $h^{\text{th}}$
dimension, and $\bm{N}_{\ell}^{h, s}\coloneqq({N}_{\ell}^{1,
  s},...,{N}_{\ell}^{D, s}),\bm{N}_{\ell}^{h, m}$, and
$\bm{N}_{\ell}^{h, c}$ are the $D$-dimensional expansion orders right
before the $\ell^{\text{th}}$ scaling, moving, or change of expansion
order in the $h^{\text{th}}$ dimension. ${\b'}^{h, s}_{\ell}$ is the
ratio $\tilde{\b}^{h,s}_{\ell}/\b^{h, s}_{\ell}$ where $\tilde{\b}^{h,
  s}_{\ell}$ is the scaling factor after the $\ell^{\text{th}}$
scaling in the $h^{\text{th}}$ dimension, $d_{\ell}^{h, m}\coloneqq
|{{\tilde{x_0}}}_{\ell}^{h, m} - {x_0}_{\ell}^{h, m}|$
(${\tilde{x_0}}_{\ell}^{h, m}$ is the new displacement) is the
absolute value of the change in displacement in the $\ell^{\text{th}}$
moving step in the $h^{\text{th}}$ dimension, and
$\tilde{N}_{\ell}^{h, c}$ is the expansion order after the $\ell^{th}$
changing the expansion order in the $h^{\text{th}}$ dimension. $t_i$
is the time for carrying out the $i^{\text{th}}$ scaling, moving, or
$p$-adaptive technique in any dimension and if within the same time
step more than one of those techniques in any dimension is used, those
$t_i$ may be the same but are listed in the order of carrying out
those techniques.

Equation~\eqref{multidimension} can be proved in a
dimension-by-dimension manner to evaluate the error caused by scaling
Eq.~\eqref{scaleerror}, moving Eq.~\eqref{moving}, and coarsening
Eq.~\eqref{refinement}. As with Eq.~\eqref{errorbound}, we also
conclude that in multi-dimension case the optimal strategy for
choosing parameters is to set $q^h \lesssim 1$ and $0 \lesssim
\delta^h$ in each dimension so that the change in the scaling factor
or the displacement results in numerical accuracy but does not result
in over-scaling or over-shifting. From the error lower bound in
Eq.~\eqref{lowerbound}, $1 \lesssim \nu^h$ and $1 \lesssim \mu^h$ are
required so that $\b^h$ and $x_0^h$ are adjusted in each dimension $h$
without incurring too large a projection error.

As for coarsening across higher dimensions, a larger $\eta_0^{h}$
could lead to a larger minimal expansion order in each dimension and
improve accuracy, but larger expansion orders lead to higher
computational cost, especially for high-dimensional problems (as the
total number of coefficients are $\Pi_{h=1}^D N^h$). Similarly,
decreasing the initial refinement threshold $\eta^h$ or $\gamma^h$, or
the adjustment ratio $\eta^h$ in the $h^{\text{th}}$ direction, will
lead to smaller errors and higher computational costs.

\subsection{Prior error estimate}

In addition to the posterior upper error bound of
Eq.~\eqref{errorbound}, we can also derive a prior error upper bound
of using the adaptive spectral method to solve Eq.~\eqref{modelpro} in
which the error estimate only depends on the solution itself.  First,
for the scaling technique, when we change the scaling factor from $\b$
to $\tilde{\b}$ and use $\pi_{N, x_0}^{\tilde{\beta}}U_{N,
  x_0}^{\beta}$ as the new numerical solution, the error $\|u(\cdot,
t)-\pi_{N, x_0}^{\tilde{\beta}}U_{N, x_0}^{\beta}(\cdot, t)\|$ can be
bounded by

\begin{equation}
\begin{aligned}
\|u(\cdot, t)-\pi_{N, x_0}^{\tilde{\beta}}U_{N, x_0}^{\beta}(\cdot, t)\| 
& \leq \|(I-\pi_{N, x_0}^{\tilde{\beta}})u(\cdot, t)\| 
+ \|\pi_{N, x_0}^{\tilde{\beta}}(u - U_{N, x_0}^{\beta})(\cdot, t)\|,\\
\: & \leq \|(I-\pi_{N, x_0}^{\tilde{\beta}})u(\cdot, t)\|+\|u(\cdot, t) 
- U_{N, x_0}^{\beta}(\cdot, t)\|.
\end{aligned}
\label{scaling1}
\end{equation}
In Eq.~\eqref{scaling1}, the term $\|(I-\pi_{N,
  x_0}^{\tilde{\beta}})u(\cdot, t)\|$ in the last equation is the
increment in the error bound resulting from scaling (scaling
error). Similarly, if we carry out the moving technique and change the
displacement of the basis function from $x_0$ to $\tilde{x}_0$ and use
$\pi_{N, \tilde{x}_0}^{\beta}U_{N, x_0}^{\beta}$ as the new numerical
solution, the error $\|u-\pi_{N, \tilde{x}_0}^{\beta}U_{N,
  x_0}^{\beta}\|$ can be bounded by

\begin{equation}
\begin{aligned}
\|u(\cdot, t)-\pi_{N, \tilde{x}_0}^{\beta}U_{N, x_0}^{\beta}(\cdot, t)\|  & 
\leq \|(I-\pi_{N, \tilde{x}_0}^{\beta})u(\cdot, t)\| 
+ \|\pi_{N, \tilde{x}_0}^{\beta}(u - U_{N, x_0}^{\beta})(\cdot, t)\| \\
\: & \leq \|(I-\pi_{N, \tilde{x}_0}^{\beta})u(\cdot, t)\|
+\|u - U_{N, x_0}^{\beta}(\cdot, t)\|.
\end{aligned}
\label{moving2}
\end{equation}
As for the $p$-adaptive technique, refinement will not bring any
additional error since $\pi_{\tilde{N}, x_0}^{\b}U_{N, x_0}^{\b} =
U_{N, x_0}^{\b}, \tilde{N}>N$. However, the error after coarsening and
using $\pi_{\tilde{N}, x_0}^{\b}U_{N, x_0}^{\b},\tilde{N}<N$ to
replace the original numerical solution $U_{N, x_0}^{\b}$ can be
bounded by

\begin{equation}
\begin{aligned}
\|u(\cdot, t) - \pi_{\tilde{N}, x_0}^{\b}U_{N, x_0}^{\b}(\cdot, t)\| 
& \leq\|u(\cdot, t)-\hat{U}_{N, x_0}^{\b}(\cdot, t)\|
+\| (\pi_{N, x_0}^{\b} - \pi_{\tilde{N}, x_0}^{\b})u(\cdot, t)\| \\
\: & \leq \|u(\cdot, t) - U_{N, x_0}^{\beta}(\cdot, t)\| 
+ \|(\pi_{N, x_0}^{\b}-\pi_{\tilde{N}, x_0}^{\b})u(\cdot, t)\|
\end{aligned}
\label{refinement2}
\end{equation}
where 

\begin{equation}
\hat{U}_{N, x_0}^{\b} = \pi_{\tilde{N}, x_0}^{\b}U_{N, x_0}^{\b} 
+ \sum_{i=\tilde{N}+1}^{N}\hat{u}_{i, x_0}^{\b} 
\hat{\mathcal{H}}_{i, x_0}^{\b}(x),\,\,\, \hat{u}_{i, x_0}^{\b} 
= (u(\cdot, t), \hat{\mathcal{H}}_{i, x_0}^{\b}(\cdot)).
\end{equation}

Finally, as with the derivation of Eq.~\eqref{errorbound}, we can
obtain an error bound which only depends on the solution $u$

\begin{equation}
\begin{aligned}
e(T)\leq & \sum_{i=0}^{K^s+K^m+K^c} \frac{b_{N_i, \b_i}+B_0}{b_{N_i, \b_i}}
\|(I-\pi_{N_i, {x_0}_i}^{\beta_i})u\|_{X(t_i, t_{i+1})}  \\
\: & \qquad +\sum_{\ell=1}^{K^s}\|(I- 
\pi_{N_{\ell}^s, x_{L_{\ell}^s}}^{\tilde{\beta}_{\ell}^s})u(\cdot, t_{\ell}^s)\|\\
\: & \qquad + \sum_{\ell=1}^{K^m}\|(I - 
\pi_{N_q^m, \tilde{x_0}_{\ell}^m}^{\beta_{\ell}^m})u(\cdot, t_{\ell}^m)\| \\
\: & \qquad + \sum_{\ell=1}^{K^c}\|(\pi_{N, x_0}^{\b}
- \pi_{\tilde{N}_{\ell}^c, {x_0}_{\ell}^c}^{\beta_{\ell}^c})u(\cdot, t_{\ell}^c)\|.
\end{aligned}
\label{errorbound2}
\end{equation}
Therefore, the posterior error estimate Eq.~\eqref{errorbound} gives
us more information on how we should choose the parameters in the
adaptive techniques to determine $\beta, x_0, N$. Prior error
bounds for adaptive spectral methods for $(D+1)$-dimensional model
problems ($x\in\mathbb{R}^D$) can be straightforwardly derived which takes a similar form of
Eq.~\eqref{errorbound2} but is excluded for brevity.

{\color{blue}\subsection{Frequency indicator and lower error bound}}

\added{As proposed in \cite{xia2020b,xia2020a}, the major goal of
  implementing our adaptive techniques is to maintain a small
  frequency indicator as defined in Eq.~\eqref{f_indicator}.  Here, we
  explicitly show that the frequency indicator is closely
  related to the error and why controlling it leads to accurate
  implementation of our adaptive techniques.  Actually, from
  Eq.~\eqref{f_indicator} we have}

\begin{equation}
\mathcal{F}(U_{N, x_0}^{\b})(\|u(\cdot, t)\| - e(t)) 
\leq \|(I- \pi_{N-M, x_0}^{\b}) u(\cdot, t)  \| + e(t), 
\end{equation}
which implies

\begin{equation}
\begin{aligned}
e(t)&\geq \frac{\mathcal{F}(U_{N, x_0}^{\b}(x, t))\|u(\cdot, t)\| 
- \|(I- \pi_{N-M, x_0}^{\b}) u(\cdot, t)  \|}{1 + \mathcal{F}(U_{N, x_0}^{\b})}\\
&\approx \mathcal{F}(U_{N, x_0}^{\b})\|u(\cdot, t)\| 
- \|(I- \pi_{N-M, x_0}^{\b}) u(\cdot, t)\|
\end{aligned}
\label{e_freq}
\end{equation}
when the frequency indicator $\mathcal{F}(U_{N, x_0}^{\b})=o(1)$ for
any $t$. \added{Therefore, the relationship between the lower error bound
and the frequency indicator is nearly linear, and thus monitoring and
controlling it leads to a small lower error bound.}

Since a function that decays more slowly or is more oscillatory as
time increases tends to have a larger frequency indicator, as shown in
\cite{xia2020b,xia2020a}, one should dynamically
%
%
switch to basis functions that decay more slowly, or
%
%
incorporate more oscillatory basis functions. Therefore, in the
adaptive spectral method shown in Fig.~\ref{algmovingscaling},
controlling the frequency indicator is achieved by either decreasing
the scaling factor (``Scale'') or increasing the expansion order
(``Refine''). 

\added{In the scaling and $p$-adaptive techniques, the scaling
  threshold $\nu$ for the scaling technique, the initial threshold
  $\eta_0$ for refining, as well as the ratio of the post-refinement
  adjustment factor $\gamma$ defined in
  Table~\ref{tab:model_variables} determine the tolerable rate of
  increase in the frequency indicator between two consecutive
  timesteps. Thus, we again justify that setting $\nu \gtrsim1$,
  $\eta\gtrsim 1$, and $\gamma\gtrsim 1$ can suppress increases in the
  frequency indicator, thus effectively suppressing the lower error
  bound if $\|u(\cdot, t)\|$ is uniformly bounded for $t\in[0,
    T]$. Because Eq.~\eqref{e_freq} does not depend on the underlying
  model or the numerical discretization, controlling the frequency
  indicator works well within adaptive spectral methods applied in a
  variety of different models.}

\added{On the other hand, as the error tends to accumulate over time,
  it is usually the case that}

\begin{equation}
e(T)\gtrsim\max_{0\leq t\leq T} \mathcal{F}(U_{N, x_0}^{\b}(x, t))
\|u(\cdot, t)\| - \|(I- \pi_{N-M, x_0}^{\b}) u(\cdot, t) \|.
\end{equation}
\added{Therefore if the frequency indicator decreases, one can
  consider increasing the scaling factor or reducing the number of
  basis functions allowing for modest increases in the frequency
  indicator. As long as the frequency indicator does not surpass the
  frequency indicator in previous timesteps, the error bound remains
  unchanged under the assumption that $\|u(\cdot, t)\|$ does not
  change significantly over time. By increasing the scaling factor,
  allocation points are more densely distributed making it possible to
  reduce their number via coarsening and to improve computational
  efficiency by using fewer basis functions.}

\section{Numerical results}
\label{Numericalex}
In our numerical examples, we numerically solve Eq.~\eqref{modelpro}
by discretizing time according to $t_j=j\Delta{t}$ and using the
scheme Eq.~\eqref{time_forward} to forward time from $t_j$ to
$t_{j+1}$.  Adaptive techniques will be used to adjust the basis
functions at different timesteps $t_j$.  The matrix-vector product
$e^{-{\bm A}_N^{\beta}(t_{j+1}-t_j)}\bm{u}_{N, x_0}^{\b}(t_{j})$ in
Eq.~\eqref{time_forward} is calculated using a ``scaling and squaring''
method in \citep{1978Nineteen}, \textit{i.e.}, we rewrite

\begin{equation}
e^{-{\bm A}_N^{\beta}(t_{j+1}-t_j)}\bm{u}_{N, x_0}^{\b}(t_j) 
= \big(e^{-\tfrac{{\bm A}_N^{\beta}(t_{j+1}-t_j)}{3}}\big)^3\bm{u}_{N, x_0}^{\b}(t_j)
\end{equation}
and evaluate $e^{-\tfrac{{\bm A}_N^{\beta}(t_{j+1}-t_j)}{3}}\bm{u}_{N,
  x_0}^{\b}(t_j)$ by Taylor expansion. The integral
$\int_{t_j}^{t_{j+1}}\!\! e^{-{\bm
    A}_N^{\beta}(t_{j+1}-t_j)}\bm{F}_{N, x_0}^{\b}(t)\dd{t}$ on the
RHS of Eq.~\eqref{time_forward} is evaluated by the Gauss-Legendre
formula described in \citep{xia2020a}.
In all examples, the error denotes the relative $L^2$-error
\begin{equation}
\frac{\|u(\cdot, t) - U_{N, x_0}^{\beta}(\cdot, t)\|}{\|u(\cdot, t)\|}.
\end{equation}
First, we numerically investigate how the parameters of the
  scaling and moving techniques affect the performance of the adaptive
  spectral method and numerically verify the conclusions drawn from
  Eq.~\eqref{errorbound}, namely, to set $q\lesssim 1, 1\lesssim\nu$
  for scaling, and $0\lesssim\delta, 1\lesssim\mu$ for moving in order
  to accurately adjust the scaling factor and translation of the basis
  functions. We also wish to explore how to appropriately set the
  parameters in the $p$-adaptive technique, the refinement threshold
  $\eta$, the coarsening threshold $\eta_0$, and the $\eta$ adjustment
  ratio to achieve higher accuracy while reducing the computational
  cost. In this work, all computations were performed using Matlab
  R2017a on a laptop with a 4-core Intel(R) Core(TM) i7-8550U CPU @
  1.80 GHz.

\begin{example}
\rm

We consider solving the following parabolic equation in the weak form
\begin{equation}
\begin{aligned}
& (u_t(x, t), v) + \big(u_{x}(x, t), v_x(x, t)\big)  = \big(f(x, t), v(x, t)\big),\,\, 
\forall v(x)\in H^1(\mathbb{R}), \,\, u(x, 0) 
= e^{\text{i}x}e^{-x^{2}/4},\\
& \qquad  f(x, t) = \frac{(x-2t) + (t+1)^3 + 2\text{i}(x-t)(1+t)}{(t+1)^{\sfrac{3}{2}}}
\exp\left[\text{i}(t+1)x-\frac{(x - 2t)^2}{4(t+1)}\right] 
\end{aligned}
\label{PDE0}
\end{equation}
which admits an analytic solution
\begin{equation}
u(x,t) = \frac{1}{\sqrt{t+1}}\exp\left[\text{i}(t+1)x-\frac{(x - 2t)^2}{4 (t+1)}\right].
\end{equation}
Not only is the center of the solution translating rightward at speed
$2t$, its magnitude $|u(x, t)|=\frac{1}{\sqrt{t+1}}
\exp\left(-\frac{(x - 2t)^2}{4(1+t)}\right)$ decays more slowly for
larger $\vert x\vert$. The solution also incurs higher frequency
spatial variations as time increases due to the
$\exp\left(\text{i}(t+1)x\right)$ factor.
\begin{figure}[tbhp]
\begin{center}
      \includegraphics[width=5in]{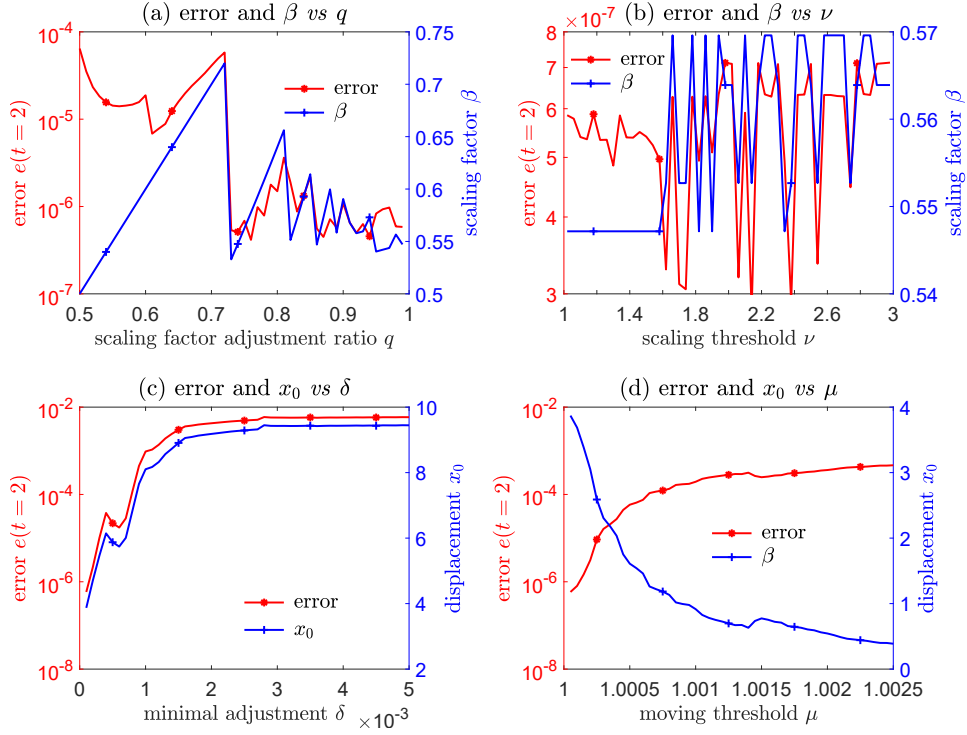}
\end{center}
        \caption{\footnotesize Plots of the error at $t=2$ and the scaling
          factor $\beta$ or the displacement $x_0$ when tuning the
          scaling factor adjustment ratio $q$ and the scaling
          threshold $\nu$ or the minimum displacement $\delta$ and the
          moving threshold $\mu$. (a) The error tends to be smaller as
          $q$ decreases to $1$, indicating that $q\lesssim 1$ is
          crucial for proper adjustment of the scaling factor. (b) As
          $\nu$ is increased, the scaling technique could be impeded,
          but the error is not very sensitive to $\nu$ if $q$ is
          small. (c) The error is strongly correlated with $x_0$ and a
          large $\delta$ can lead to over-adjustment of the
          displacement $x_0$, resulting in a larger error. (d) A large
          $\mu$ will make it harder to activate the moving technique,
          leading to a smaller $x_0$ and a larger error.}
     \label{fig0}
\end{figure}
Therefore, all three adaptive techniques are expected to be required.
Upon setting $\Delta{t} = 2\times 10^{-4}$ and solving
Eq.~\eqref{PDE0} up to $t=2$, we investigate how the parameters in the
three adaptive techniques affect performance. The initial scaling
factor, displacement, and expansion order are set to $\beta=1$,
$x_0=0$, and $N=40$. First, we test how the scaling threshold $\nu$,
the scaling factor adjustment ratio $q$, the moving $\mu$, and the
minimum displacement step $\delta$ affect the performance of the
scaling and moving techniques. We keep the expansion order fixed since
it has been illustrated that the effects of improper scaling or moving
can be offset by increasing the expansion order $N$ but at the expense
of increased computational cost \citep{xia2020b}. Initially, we set
the parameters $q=0.99, \nu=1.02, \delta=10^{-4}$, and $\mu=1.00005$,
and then change each of them one at a time. Imposing the maximal
allowable displacement within each timestep $d_{\max}=0.01$, the upper
scaling factor limit $\overline{\b}=0.2$, and lower scaling factor
limit $\overline{\b}=5$, we plot the relative $L^2$-error $e(t=2)$
along with the scaling factor when we change $q$ and $\nu$, and we
plot $e(t=2)$ along with $x_0$ when we change $\delta$ and $\mu$.

Fig.~\ref{fig0}(a) shows that $q\lesssim 1$ is required for the
scaling technique to properly adjust the scaling factor. When
$q\lesssim 1$ and we vary $\nu$ from $1$ to $2$, the error, as well as
the scaling factor $\b$, do not change much, indicating that the
scaling technique is more sensitive to $q$ than to $\nu$.  Therefore,
keeping $q\lesssim 1$ is more important than keeping $1\lesssim
\nu$. Fig.~\ref{fig0}(c) shows that the error is highly correlated
with $x_0$, suggesting that it is critical to properly move the basis
functions to capture the displacement of the solution.  Having
$0\lesssim\delta$ is important so that the displacement $x_0$ is not
over-adjusted. Finally, as shown in Fig.~\ref{fig0}(d), increasing
$\mu$ will make the moving technique less sensitive to the translation
of the basis functions and lead to a larger error. Thus,
$1\lesssim\mu$ is recommended for the moving technique.
\begin{figure}[tbhp]
\begin{center}
      \includegraphics[width=5in]{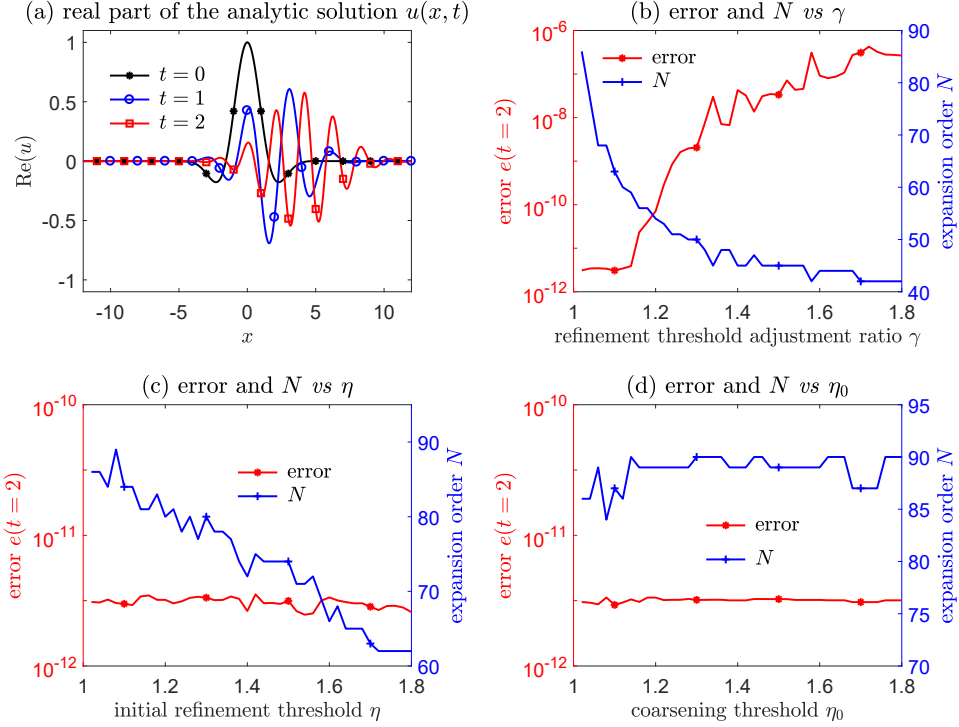}
\end{center}
\caption{\footnotesize Plots of the real part of the analytic solution
  $\text{Re}(u)(x, t)$ at different times, the error and the expansion
  order $N$ at $t=2$ when we vary the refinement threshold adjustment
  ratio $\gamma$, the initial refinement threshold $\eta$, and the
  coarsening threshold $\eta_0$. (a) The real part of the analytic
  solution, which translates rightward, becomes more diffusive, and is
  increasingly oscillatory over time. (b) The error increases with
  $\gamma$ while the expansion order decreases with $\gamma$. A larger
  $\gamma$ implies a faster-increasing refinement threshold
  $\eta$. (c) A larger initial refinement threshold $\eta$ results in
  a smaller expansion order at $t=2$, yet the error is not reduced as
  $\eta$ decreases and $N$ increases with the initial $\gamma$. This
  indicates that as long as $\gamma$ is small enough, a larger initial
  $\eta$ can be tolerated to lead to a smaller computational cost
  without compromising accuracy. (d) The expansion order $N$ tends to
  increase as the coarsening threshold $\eta_0$ increases.}
     \label{fig00}
\vspace{-0.1in}
\end{figure}
Next, we investigate how the initial refinement threshold $\eta$, the
refinement threshold adjustment ratio $\gamma$, and the coarsening
threshold $\eta_0$ affect the $p$-adaptive technique's performance
when $q=0.99, \nu=1.02, \delta=10^{-4}$, and $\mu=1.00005$ are fixed,
and the initial variables are set to $\b=1, x_0=0, N=40$.  Fixing the
maximum increment to $N_{\max}=6$, we start with the initial parameter
values $\gamma=1.02$, $\eta=1.05$, and $\eta_0=1.02$, and vary each of
them one by one and plot the relative $L^2$-error and $N$.
Fig.~\ref{fig00}(a) shows that apart from translating rightward and
decaying more slowly, the analytic solution is increasingly
oscillatory which requires adjusting the expansion order $N$ of the
numerical solution. Fig.~\ref{fig00}(b) shows that if $\gamma$ is
large, then the threshold for increasing the expansion order $\eta$
will increase more quickly. This renders the $p$-adaptive technique
unable to sufficiently adjust the expansion order, leading to smaller
expansion orders $N$ and larger errors.  Fig.~\ref{fig00}(c) shows
that the larger the initial threshold $\eta$ for increasing the
expansion order, the smaller the expansion order.  In the depicted
regime, larger initial values of $\eta$ do not degrade accuracy since
$N\gtrsim 65$ is sufficient to maintain high accuracy. Therefore, to
maintain accuracy while reducing the computational burden, it is
crucial to set $1\lesssim \gamma$ so that the $p$-adaptive technique
can capture oscillatory behavior over long periods of time. Using a
smaller initial $\eta$ may lead to more computational costs but does
not lead to improvement in accuracy. Overall, since the function
exhibits higher frequency spatial oscillations as time increases,
coarsening is typically not activated. However, a large coarsening
threshold $\eta_0$ can still impede coarsening, resulting in a slightly
larger $N$ than a smaller $\eta_0$ (Fig.~\ref{fig00}(d)).

Finally, as shown in Figs.~\ref{fig0} and~\ref{fig00}, we numerically
verify that the appropriate strategy for the adaptive spectral
parameters is to set $q\lesssim 1, 1\lesssim\nu, 0\lesssim\delta$, and
$1\lesssim\mu$. In fact, for good performance, the scaling procedure
strongly requires $q\lesssim 1$ and the moving procedure requires both
$0\lesssim\delta$ and $1\lesssim\mu$.  For an effective refinement, it
is more important to set $1\lesssim \gamma$ rather than to set the
initial $1\lesssim\eta$ (\textit{i.e.}, setting $1\lesssim\gamma$
rather than setting the initial $1\lesssim\eta$ leads to more accurate
results with smaller computational costs).
%
\end{example}
\begin{figure}[tbhp]
\begin{center}
      \includegraphics[width=5in]{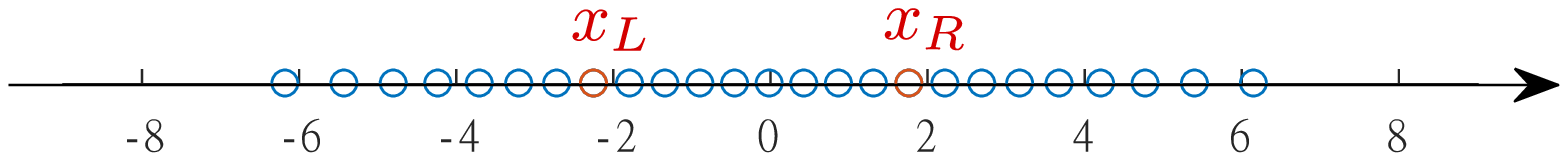}
\end{center}
        \caption{\footnotesize Distribution of the collocation points of
          generalized Hermite functions $\{\hat{\mathcal{H}}_{i,
            x_0}^{\b}\}_{i=0}^N$ with $\b=1, x_0=0$, and
          $N=24$. $x_L\coloneqq x_{[\f{N}{3}]}^{\b}$ and $x_R
          \coloneqq x_{[\f{2N+2}{3}]}^{\b}$ are marked in red. The
          number of collocation points that are in the right-exterior
          region $(x_R, \infty)$ for calculating $\mathcal{E}_R$ and
          in the left-exterior region $(-\infty, x_L)$ for calculating
          $\mathcal{E}_L$ are both approximately $N/3$.}
     \label{fig_hermite}
\end{figure}
When using the generalized Hermite functions defined in $\mathbb{R}$,
the desired solution might move leftward or rightward, requiring both
leftward and rightward displacement of the basis functions. Since only
rightward basis function shifts have been previously considered
\citep{xia2020b,xia2020a}, here, we generalize the moving technique to
allow for bidirectional adjustment of the displacement $x_0$.
\added{It has been proposed that controlling an exterior-error
  indicator leads to small errors in the exterior domain, relative to
  the total error, resulting in better approximation of the solution
  in the exterior region.  Therefore, bidirectional moving might
  maintain relatively small errors in both left- and right-exterior
  regions of $\mathbb{R}$.} We first propose a left exterior-error
indicator
 \begin{equation}
   \mathcal{E}_L(U_{N,x_0}^{\b}) = \f{\|\p_x U_{N,x_0}^{\b}
\cdot\mathbb{I}_{(-\infty,x_L)}\|}{\|\p_x U_{N,x_0}^{\b}\cdot
     \mathbb{I}_{(-\infty,+\infty)}\|},
\label{errorindicator}
\end{equation}
where we use $x_L = x_{[\f{N}{3}]}^{\b}$ following the often-used
$\tfrac{2}{3}$-rule \citep{Hou2007Computing,Orszag1971On}. The left
exterior-error indicator~\eqref{errorindicator} can be seen as the
upper bound for the ratio of the error in $(-\infty, x_L)$ to the
error across the whole space $\mathbb{R}$, in analogy to the (right)
exterior-error indicator $\mathcal{E}(U_{N,x_0}^{\b})$ defined in
Eq.~\eqref{errorindicator2}, which we shall denote below by
$\mathcal{E}_R(U_{N, x_0}^{\b})$ for clarity.
%
%
The number of nodes in the left-exterior region $(-\infty, x_L)$ and
in the right-exterior region $(x_R, \infty)$ are both roughly
$\frac{N}{3}$.  It was shown in \citep{xia2020a} that if the right
exterior-error indicator~\eqref{errorindicator2} increases, then the
ratio of the error in the right exterior region $(x_R, +\infty)$ to
the total error may also increase, suggesting that one should move the
basis functions rightward (increase $x_0$). In Fig.~\ref{fig_hermite},
we show the positions of collocation nodes of generalized Hermite
functions $\{\mathcal{H}_{i, x_0}^{\b}\}_{i=0}^N$ with $\b=1, x_0=0$,
and $N=24$. The endpoints $x_L$ and $x_R$ are shown in red, showing
that the right and left exterior regions, $(x_R, \infty)$ and $(
-\infty, x_L)$, are near-symmetric. The left exterior-error indicator
~\eqref{errorindicator} also measures the ratio of the error in the left
exterior region $(-\infty, x_L)$ to the total error, and, if it 
increases, one can consider shifting the basis functions leftward (decrease $x_0$). With
both left and right exterior-error indicators, we propose the
following bidirectional moving scheme.

\begin{algorithm}[t]
  \caption{\small Pseudo-code of the bidirectional
    exterior-error-dependent moving technique.}
\begin{algorithmic}[1]
\State Initialize $N$, $\Delta t$, $T$,  $\beta$, $x_0$, $U_{N,x_0}^{\beta}(x, 0)$, $\mu>1$, $d_{\rm max}>\delta>0$
\State $t\gets 0$
\State $x_R \gets x^{\beta}_{[\f{2N+2}{3}]}$
\State $x_L \gets x^{\beta}_{[\f{N}{3}]}$
\State $\tilde{\mathcal{E}}_R \gets \Call{right\_exterior\_error\_indicator}{U_{N,x_0}^{\b}(x, 0)}$\label{m:e0}
\State $\tilde{\mathcal{E}}_L \gets \Call{left\_exterior\_error\_indicator}{U_{N,x_0}^{\b}(x, 0)}$\label{m:e1}
\While{$t<T$}
\State $U_{N,x_0}^{\beta}(x, t+\Delta t)\gets \Call{evolve}{U_{N,x_0}^{\beta}(x, t),\Delta t}$
\State $\mathcal{E}_R\gets \Call{right\_exterior\_error\_indicator}{{U_{N,x_0}^{\beta}(x, t+\Delta{t})}}$\label{m:e2}
\State $\mathcal{E}_L\gets \Call{left\_exterior\_error\_indicator}{{U_{N,x_0}^{\beta}(x, t+\Delta{t})}}$\label{m:e3}
\If{$\mathcal{E}_R > \mu \tilde{\mathcal{E}}_R ~|| ~\mathcal{E}_L>\mu \tilde{\mathcal{E}}_L$}\label{m:if}
\State $d_R\gets \Call{move\_right}{U_{N,x_0}^{\beta}(t+\Delta t), \delta, d_{\rm max}, \mu e_0}$\label{m:d0}
\State $d_L\gets \Call{move\_left}{U_{N,x_0}^{\beta}(t+\Delta t), \delta, d_{\rm max}, \mu e_1}$\label{m:d1}
\State $U_{N, x_0}^{\b}(x, t)\gets \pi_{N, x_0+d_R-d_L}^{\b}U_{N, x_0}^{\b}(x, t+\Delta t)$
\State $x_0 \gets x_0+d_R-d_L$\label{m:x0}
\State $x_L \gets x_L+d_R-d_L$\label{m:xL}
\State $x_R \gets x_R+d_R-d_L$\label{m:xR}
\State $\tilde{\mathcal{E}}_R \gets \Call{right\_exterior\_error\_indicator}{{U_{N,x_0}^{\beta}(x, t+\Delta{t})}}$\label{m:e4}
\State $\tilde{\mathcal{E}}_L \gets \Call{left\_exterior\_error\_indicator}{{U_{N,x_0}^{\beta}(x, t+\Delta{t})}}$\label{m:e5}
\EndIf
\State $t\gets t+\Delta t$
\EndWhile
\end{algorithmic}\label{algmoving}
\end{algorithm}
In Alg.~\ref{algmoving}, the
$\Call{left\_exterior\_error\_indicator}{}$ subroutine calculates the
left exterior-error indicator by Eq.~\eqref{errorindicator} and the
$\Call{right\_exterior\_error\_indicator}{}$ calculates the right
exterior-error indicator by Eq.~\eqref{errorindicator2}. If the right
or left exterior-error indicator is larger than their corresponding
thresholds, $\textit{i.e}$, $\mathcal{E}_R>\mu \tilde{\mathcal{E}}_R$
or $\mathcal{E}_L>\mu \tilde{\mathcal{E}}_L$, the moving technique is
activated, calculating the rightward displacement $d_0$ or the
leftward displacement $d_1$ of the basis functions. In
\citep{xia2020a}, the rightward displacement $d_R=\min\{n_R\delta,
d_{\rm max}\}$ is determined by the $\Call{move\_right}{}$ subroutine
in Line~\ref{m:d0}, where $n$ is the smallest integer satisfying
$\mathcal{E}_R(U_{N,x_0}^{(\a, \b)}(x-n_R\delta, t)) <\mu
\tilde{\mathcal{E}}_R$. Similarly, the leftward displacement $d_L =
\min\{n_L\delta, d_{\rm max}\}$ is determined by the
$\Call{move\_left}{}$ subroutine in Line~\ref{m:d1}, where $n_L$ is
the smallest integer satisfying $\mathcal{E}_L(U_{N,x_0}^{(\a,
  \b)}(x+n_L\delta, t)) <\mu \tilde{\mathcal{E}}_L$. \added{Notice
  that the error estimate of the adaptive spectral method in
  Theorem~\ref{theorem1} does not depend on the direction of
  displacements. Therefore it applies to both the bidirectional moving
  technique Alg.~\ref{algmoving} and the one-sided moving technique
  proposed in \citep{xia2020b}.}

\begin{example}
\rm Consider numerically solving the following parabolic equation in the weak form in
$\mathbb{R}\times\mathbb{R}^+$

\begin{equation}
\begin{aligned}
(u_t(x, t), v) + \big(u_{x}(x, t), v_x(x, t)\big)  = \big(f(x, t), v(x, t)\big),\,\,\, 
\forall v(x)\in H^1(\mathbb{R}), \,\, u(x, 0) = e^{-x^{2}} \sin x,
\end{aligned}
\label{PDE1}
\end{equation}
where 
\begin{equation}
f(x, t) = \left\{\begin{array}{ll}
\begin{array}{l}\bigg[\Big(3-2(x+vt)\big(v+2(x+vt)\big)\Big)\sin(x+vt)\\
\hspace{1.6cm}+\Big(v+4(x+vt)\Big)\cos(x+vt)\bigg]
e^{-(x+vt)^{2}}\end{array} & t\leq 2, \\[22pt]
\begin{array}{l}
\bigg[\Big(3 - 4\big(x+v(4-t)\big)^2+2v\big(x + v(4-t)\big)\Big)\sin(x-v(t-4))\\
\hspace{1.2cm} +\Big(4x+v(15-4t)\Big)\cos(x-v(t-4))\bigg]
e^{-(x-v(t-4))^{2}}\end{array} & t\geq 2. \end{array}\right.
\label{F1}
\end{equation}
This PDE is solved by

\begin{equation}
u(x, t) =\left\{\begin{array}{ll}
\displaystyle  e^{-(x+vt)^2}\sin\left(x+vt\right)  & t\leq2 ,\\[13pt]
\displaystyle  e^{-(x-vt+4v)^2}\sin\left(x-vt+4v\right) &  t \geq 2.
\end{array}\right.
\end{equation}

We set $v=2$ in Eq.~\eqref{F1} so that the center of the solution
moves with velocity $-2$ from $x=0$ to $x=-4$ when $t\in[0, 2]$, and
when $t\in[2, 6]$ the center of the solution moves from $x=-4$ to
$x=4$ with velocity $+2$. Since the solution displays only convective
behavior, we deactivate the scaling and $p$-adaptive procedures and
apply only the moving technique.  Since the translation switches from
leftward to rightward at $t=2$, the moving technique needs to allow
for both leftward and rightward displacement of the basis functions.
The parameters in the moving technique are set to be $\mu=1.0005,
\delta=0.0005$, and the maximal displacement within a timestep
$d_{\max}=0.2$.  We take the scaling factor, the expansion order, and
the initial displacement of the basis function to be $\b_0=1.2,
N_0=24, x_0=0$, respectively, and plot the results obtained with no
moving technique, the leftward-only moving technique, the
rightward-only moving technique, and the bidirectional moving
technique.

\begin{figure}[tbhp]
\begin{center}
      \includegraphics[width=5in]{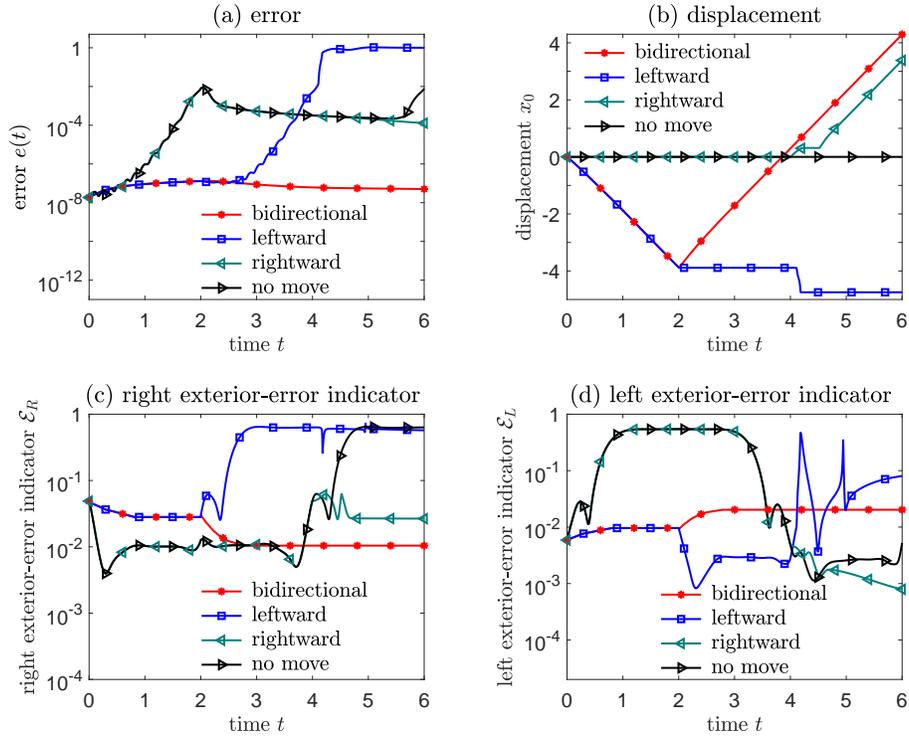}
\end{center}
\caption{\footnotesize Plots of the error, $x_0$, the left
  exterior-error indicator Eq.~\eqref{errorindicator}, and the right
  exterior-error indicator Eq.~\eqref{errorindicator2}. (a) The
  bidirectional moving technique Alg.~\ref{algmoving} can main the
  smallest error while failure to accommodate either leftward or
  rightward displacement leads to much larger errors. (b,c,d) The
  displacement $x_0$, the left exterior-error indicator, and the right
  exterior-error indicator of spectral methods with the bidirectional,
  the leftward-only, the rightward-only moving technique, and the
  spectral method without any moving.}
\label{fig2}
\end{figure}
Fig.~\ref{fig2}(a) shows that the spectral method equipped with the
bidirectional moving technique (red) can maintain the smallest error
because the displacement $x_0$ can be decreased when $t\in[0, 2]$ and
increased when $t>2$ (see Fig.~\ref{fig2}(b)). The spectral method
with the leftward-only moving technique (blue) can maintain a small
error in $[0, 2]$ when the center of the function moves leftward but
fails to keep the error small when $t>2$ due to its inability to
increase $x_0$. When $t<2$, the rightward-only moving technique
(green) cannot decrease the displacement $x_0$ and therefore the error
for the rightward-only moving technique is large at
$t=2$. Furthermore, large error accumulation before $t=4$ of the
rightward-only moving technique makes it unable to properly increase
$x_0$ for $t>4$ when the center of the solution moves to the right of
the origin $x=0$. The right and left exterior-error indicators for the
bidirectional moving technique Alg.~\ref{algmoving} can be well
controlled as shown in Fig.~\ref{fig2}(c,d), while for the
leftward-only moving technique the right exterior-error indicator
grows dramatically when $t>2$ and for the rightward-only moving
technique, the left exterior-error indicator grows when $t<2$.
Therefore, the leftward- and rightward-only moving techniques both
fail to maintain a small error in at least one exterior region $(x_R,
\infty)$ or $(-\infty, x_L)$. The left exterior-error indicator grows
when $t<2$ (the center moves to the left of the origin) and the right
exterior-error indicator grows when $t>4$ (the center moves to the
right of the origin) for the spectral method without the moving
technique (black), suggesting that it cannot maintain a small error in
both exterior regions.

\end{example}

\section{Discussion and Conclusions}
In this paper, we carried out a numerical analysis of recently
proposed adaptive spectral methods in unbounded domains using
generalized Hermite functions.  Specifically, our analysis helps guide
parameter choice across three adaptive spectral techniques,
\textit{i.e.}, the scaling procedure, the moving procedure, and the
$p$-adaptive technique to properly adjust the three key variables
associated with these techniques, the scaling factor, the
displacement, and the spectral expansion order. Based on our analyses,
rules for properly choosing parameters in the scaling, moving, and
$p$-adaptive techniques to most efficiently and accurately solve PDEs
are derived. \added{We also explicitly explain why controlling the
  frequency indicator by using adaptive spectral methods effectively
  controls the error.} Numerical experiments were carried out to
verify our theoretical results. Furthermore, we developed a new
bidirectional moving technique to accommodate both leftward and
rightward displacements.

Even though our analysis focused on a simple parabolic model, it
nonetheless represents a first step towards understanding how adaptive
spectral methods work in solving unbounded-domain problems. \added{In
  fact, for our parabolic model, the total upper error bound is simply
  the sum of the errors from numerical discretization and from
  implementation of the adaptive schemes, providing a clear overall
  picture of errors under our adaptive spectral algorithm. Additionally,
  the lower error estimate Eq.~\eqref{e_freq} holds regardless of the
  underlying model and numerical discretization, suggesting that
  controlling a small frequency indicator always leads to a small
  lower error bound when applying adaptive spectral methods to any
  model.}

Since adaptive spectral methods have been successfully applied to
nonlinear PDEs or models containing nonlocal terms
\citep{xia2020b,xia2020a}, further analysis to explain why adaptive
spectral methods work well in these more complicated models,
particularly in unbounded domains, will be the subject of future
investigation. Understanding how adaptive spectral methods work in
complex unbounded-domain problems that arise across many disciplines
and that are computationally challenging will pave the way for their
accurate solution.

Finally, one should also perform analyses of adaptive spectral
techniques using other classes of basis functions of recent interest
\citep{tang2020rational}. These include generalized Laguerre functions
in $\mathbb{R}^+$ and the modified mapped Gegenbauer functions in
$\mathbb{R}$. Another potentially useful extension is to explore
developing methods to automatically determine and adjust the decay
rate of solutions at infinity by adaptively switching among different
classes of basis functions in order to match underlying physics or
observations.

\section*{Funding}
TC and MX were supported from the US National Science Foundation
through grant DMS-1814364. SS was supported by the National Key R\&D
program of China (No.~2020AAA0105200) and Beijing Academy of
Artificial Intelligence (BAAI).




\bibliographystyle{elsarticle-num}

\end{document}